\documentclass[12pt]{article}
\usepackage{ejorsj-s2}
\usepackage{amsmath,amssymb,amsfonts,amsthm}
\usepackage{graphicx}
\usepackage{hyperref}

\eqswitchon 
\newtheorem{theorem}{Theorem}[section]
\newtheorem{lemma}{Lemma}[section]
\newtheorem{definition}{Definition}[section]

\newcommand{\rs}[1]{{\mbox{\scriptsize \sc #1}}}
\newcommand{\vc}[1]{{\boldsymbol #1}}

\newcommand{\sr}[1]{{\cal #1}}
\newcommand{\dd}[1]{\mathbb{#1}}

\newcommand{\cp}{{\Gamma}}

\newcommand{\diag}{{\mbox{diag}}}

\newcommand{\eqn}[1]{(\ref{eqn:#1})}
\newcommand{\lem}[1]{Lemma~\ref{lem:#1}}
\newcommand{\cor}[1]{Corollary~\ref{cor:#1}}
\newcommand{\thr}[1]{Theorem~\ref{thr:#1}}

\newcommand{\rem}[1]{Remark~\ref{rem:#1}}

\newcommand{\app}[1]{Appendix~\ref{app:#1}}
\newcommand{\sectn}[1]{Section~\ref{sect:#1}}

\newcommand{\lemt}[1]{\ref{lem:#1}}

\newcommand{\thrt}[1]{\ref{thr:#1}}

\newcommand{\sect}[1]{\ref{sect:#1}}

\newcommand{\br}[1]{\langle #1 \rangle}

\newcommand{\ol}{\overline}
\newcommand{\ul}{\underline}
\newcommand{\pend}{\hfill \thicklines \framebox(6.6,6.6)[l]{}}

\newenvironment{proof*}[1]{\noindent {\sc  #1} \rm}{\pend}

\newtheorem{remark}{Remark}[section]

\newtheorem{corollary}{Corollary}[section]

\newcommand{\setsection}[2] {
\setcounter{section}{#1}
\setcounter{subsection}{0}
\setcounter{equation}{0}
\setcounter{conjecture}{0}
\setcounter{assumption}{0}
\setcounter{question}{0}
\setcounter{definition}{0}
\setcounter{theorem}{0}
\setcounter{corollary}{0}
\setcounter{lemma}{0}
\setcounter{proposition}{0}
\setcounter{remark}{0}
\setcounter{appen}{0}
\setsection*{\large \bf \thesection. #2}}

\newenvironment{mylist}[1]{\begin{list}{}
{\setlength{\itemindent}{#1mm}}
{\setlength{\itemsep}{0ex plus 0.2ex}}
{\setlength{\parsep}{0.5ex plus 0.2ex}}
{\setlength{\labelwidth}{10mm}}
}{\end{list}}

\title{DIFFUSION APPROXIMATION FOR STATIONARY ANALYSIS OF QUEUES AND THEIR NETWORKS: A REVIEW}
\author{
\begin{tabular}[h]{ccc}  
Masakiyo Miyazawa                        \\ % author names
Tokyo University of Science\\ % affiliations
\end{tabular}
}
\date{(Received October 2, 2014; Revised January 18, 2015)}
\begin{document}
%------------------------------------------------------ 
\maketitle
%-------- abstract (100 to 200 words) -------- 
\begin{abstract}
Diffusion processes have been widely used for approximations in the queueing theory. There are different types of diffusion approximations. Among them, we are interested in those obtained through limits of a sequence of models which describe queueing networks. Such a limit is typically obtained by the weak convergence of either stochastic processes or stationary distributions. We already have nice reviews and text books for them. However, this area is still actively studied, and it seems getting hard to have a comprehensive overview because mathematical results are highly technical. We try to fill this gap presenting technical background. Although those diffusion approximations have been well developed, there remains a big problem, which is difficulty to get useful information from the limiting diffusion processes. Their state spaces are multidimensional, whose dimension corresponds to the number of nodes for a single-class case and the number of customer types for a multi-class case. We now have a better view for the two dimensional case, but still know very little about the higher dimensional case. This intractability is somehow against a spirit of diffusion approximation. This motivates us to reconsider diffusion approximation from scratch. For this, we highlight the stationary distributions, and make clear a mechanism to produce diffusion approximations.
\end{abstract}
%-------- keywords (2 to 6 keywords) ------- 
\keyword{Queue, Applied probability, Generalized Jackson network, Single-class, Multi-class, Heavy traffic, Kingman's approximation, Diffusion approximation, Semi-martingale reflecting Brownian motion, Stationary distribution, Tail asymptotics}
%------------------------------------------------------ 
\section{Introduction}
\label{sect:introduction}
%------------------------------------------------------ 

The aim of this paper is to review some technical aspects on diffusion approximation arising in queues and their networks in heavy traffic. We only consider open networks which have exogenous arrivals. The diffusion approximation is meant to use a diffusion process for approximation. Here, a diffusion process is a continuous time strong Markov process which is characterized by a stochastic integral equation using a Brownian motion and continuous processes with bounded variations (e.g., see \eqn{diffusion process 1} for this integration). This diffusion process has a continuous path, and is called It\^{o} diffusion in some literature, but we will just call it a diffusion process. A comprehensive introduction about it can be found in Harrison \cite{Harr2013} (see also Kushner \cite{Kush2001} for control problems). Mathematical details can be found in Section 5 of Protter \cite{Prot2005}.

Queueing models can not be described by diffusion processes in general because customers are discrete in nature, and therefore their sample paths are not continuous.  However, analytical studies on queueing models are generally very hard particularly for queueing networks except for special cases. This motivates to use diffusion processes for their approximations because diffusion processes are characterized by fewer modeling parameters, that is, fewer primitive data, and have simpler state spaces. For example, the first two moments including covariances are typical primitive data for diffusion processes. Thus, diffusion approximation may be considered as a two-moment approximation. Their state spaces are finite dimensional vector spaces, while queueing processes require supplementary states for describing their dynamics. Those two features are analytical advantages of diffusion approximation. However, one must be careful of such simplification because important information might be lost.

Another motivation for diffusion approximation comes from queueing dynamics. In general, queues and their networks are input-output systems, and inputs and outputs can be described by additive processes. For example, a counting process of arriving customers and accumulated workload brought by customers are such additive processes. If the additive process has independent or weakly dependent increments, then it is natural to use Brownian motion for approximation. In this case, queueing characteristics such as queue lengths are functionals of Brownian motion, which may not be easily analyzed particularly for the multi-dimensional case, but fewer primitive data are advantage to study feature of the original model.

Diffusion approximation is so versatile, and the literature about it is so huge. It is hard to survey its whole area. Thus, we focus on three things. The first is about a method to derive diffusion approximations. We are only interested in those approximations supportable by theory. Namely, they should be asymptotically exact under certain limiting operations. Because objects of our interest are a stochastic process and its stationary distribution, the limiting operation is convergence in distribution under scaling of state (and time if necessary). Here, we will show that scaling in time may not be necessary if we are only concerned with the stationary distribution. This is different from the standard approach of diffusion approximation, but may not be surprising because the stationary distribution is independent of time. Furthermore, Kingman \cite{King1962} uses such scaling for the waiting time of a single server queue.

The next step is to derive approximations based on those limiting results. However, this step is often ignored in the literature because it may be considered obvious. For example, the limiting processes and distributions may be considered as approximations themselves. Here is a gap between mathematical results and approximations. To fill this gap, we go back to the original idea of Kingman \cite{King1962}. He proposes an approximation for the stationary waiting time distribution for a single server queue by moving a scaling factor into the limiting distribution. We think this is a right way to have theoretical support. Of course, we may be motivated by intuitive arguments, and should not exclude heuristics. However, it also is important to legitimize them by theory. 

The second is about a class of queueing models. Our primary interest is in queueing networks, but also consider some of single node queues because they are components of network systems. Key features of a queueing network are routing of customers and idling of servers. There are two types of routing, choosing next nodes or next types. The latter is used for a multi-class queueing network. In either case, we assume Markovian routing, which means that nodes or types are chosen with given probabilities only depending on the current nodes or types, respectively, independently of everything else. If exogenous arrivals at nodes are subject to independent renewal processes, if service times at each server are $i.i.d.$ (independent identically distributed) and if service discipline is first-come first-served, FCFS for short, then this network model is called a generalized Jackson network, GJN for short. This network is a main subject of this review.

The third is about objects for approximations. As we mentioned above, a stochastic process describing a model is a primary object. However, from the application viewpoint, its stationary distribution is also very important. Other objects are moments and tail asymptotics of the stationary distribution. Note that the weak convergence of the sequence of stochastic processes does not imply that of their stationary distributions. That is, the exchange of those limits needs proofs. If it is possible, the approximation is considered to have better support. Similarly, the convergence of the sequence of the stationary distributions does not imply the convergences of their moments and tail asymptotics. If they do so, we may say the approximation has better quality.

Theoretical support for diffusion approximation generally requires a heavy traffic condition. Under this condition, the traffic intensity at each node is close to unit, and therefore the behavior of a small queue at a fixed node is ignorable. In other words, diffusion approximation is for the behavior of a large queue. This greatly simplifies diffusion approximation for a single queue. On the other hand, in a queueing network, there are multiple queues at different nodes, and they may not simultaneously get large. For example, there may be idle servers when some queues are large. This is a crucial aspect of the queueing network, and should not be excluded in its approximation. Thus, it is interesting to see how such a behavior is implemented in diffusion approximation.

Historically, the study of diffusion approximation started from single node queues. Among them, a $GI/GI/1$ queue is known as a basic model. This model assumes that customers arrive according to a renewal process, there is a single server and service times are $i.i.d.$ Even this simple model is known to be hard to get the stationary distribution in closed form. This motivated Kingman \cite{King1962,King1965} to originate a two-moment approximation for the stationary waiting time distribution of the $GI/GI/1$ queue in heavy traffic. This approximation is generalized for a $GI/GI/s$ queue with homogeneous servers by K\"ollerstr\"om \cite{Koll1974}. Those studies directly consider the stationary distribution, so do nothing with diffusion approximation. In particular, scaling is taken only for state, as we already mentioned. However, their stationary distributions are identical with those obtained by the standard diffusion approximations (e.g., compare \eqn{W convergence 1} with \cor{GJN 1} together with \eqn{W 2}).

Independent of those studies, Borovkov \cite{Boro1965} studied diffusion approximation for the transient behavior of a many server queue. It refers to Provkov's 1963 paper, but we are unable to get this paper. However, Iglehart and Whitt \cite{IgleWhit1970a} credited that those Russian work is one of the origin of diffusion approximation. Thus, the work may be the first one to study diffusion limits for a single queue by scaling state and time. The book \cite{Boro1976} of Borovkov is a good source for those studies.

A next step was an extension from a single queue to a feedforward network. This step was firstly taken by seminal papers of Iglehart and Whitt. They obtained diffusion processes as weak limits of the sequence of suitably scaled processes, called a diffusion scaling, in \cite{IgleWhit1970,IgleWhit1970a}. Their results are basically for a single queue because functional relations of input and output of a single queue are essential. For example, they only studied a marginal process of a given node, but were unable to study a joint process for multiple nodes. However, there are two things to be notable. The first is an idea for process limit, that is, a weak limit of the sequence of stochastic processes. The second is that a reflecting Brownian motion is obtained as one of the process limits. This means that the diffusion approximation is more than a stochastic process version of the central limit theorem. This approximation is referred to as a functional limit theorem under normalization similar to that of the central limit theorem. A similar approach is also studied for a single queue in \cite{Boro1976}.

A real step toward a queueing network was taken by Harrison \cite{Harr1978}, who considered a tandem queue with two nodes, and got a two-dimensional reflecting Brownian motion on the nonnegative orthant. It still uses an explicit functional relation of input and output, but it was a big step to move from one-dimensional diffusion process to the two-dimensional one.  Harrison and Reiman \cite{HarrReim1981} studied this reflecting process for an arbitrary dimension. It is a multidimensional reflecting Brownian motion on the nonnegative orthant. Because it has a semi-martingale representation, it is called a semi-martingale reflecting Brownian motion, SRBM for short. Those studies were culminated by the work of Reiman \cite{Reim1984}. It obtains SRBM as a weak limit of suitably scaled processes for GJNs, where scaling is taken for both state and time. This scaling is called diffusion scaling, and the weak limit is called a diffusion limit.

Another source for diffusion approximation is heuristic arguments. A queueing network including GJN is analytically intractable except for some special cases. Typical examples for this exceptional case are a Jackson network and its extensions, called product form networks due to Kelly \cite{Kell1979} (see \cite{ChaoMiyaPine1999,Serf1999} for further development in this line). They are analytically tractable, but require stronger assumptions. Thus, there have been many attempts for heuristic approximations for a GJN, particularly for its stationary distribution. For example, the diffusion approximation of Kobayashi \cite{Koba1974} is one of the earliest work. Among them, Queueing Network Analyzer, QNA for short, proposed by Whitt \cite{Whit1983,Whit1983b}, is well known. This is a heuristic approximation based on the product form solution and the two-moment approximation for the queue length at each node.

The concept of QNA has been widely accepted, and similar types of approximations have been considered for more complicated network models (e.g., see \cite{Kino1997}). However, it lacks the theoretical supports which we have discussed. A theoretical improvement of QNA was studied by Harrison and Nguyen \cite{HarrNguy1990}. They propose another approximation concept, called QNET, incorporating information from a diffusion limit for the sequence of GJNs in heavy traffic. This limiting process is an SRBM. We note that this QNET is different from Q'NET (alternatively spelled as Q-NET) which has been developed for business solutions by IBM.

The SRBM is much simpler than the GJN, and its primitive data can be expressed by the mean and covariances of Brownian motion and the reflecting matrix. An approximation based on it may be considered as a two-moment approximation. However, the stationary distribution of the SRBM is still hard to get except for either one-dimensional case or the skew-symmetric case. The latter corresponds to the Jackson network, and has a similar product form solution. Thus, QNET requires approximation for the stationary distribution of the SRBM. For this, it uses either the product form or numerically approximated solution for the stationary distribution. Thus, QNET may not be fully supported by theory. Nevertheless, QNA and QNET had great impacts on theoretical studies for diffusion approximation through their design concepts and practical applications using their software packages.

In present days, diffusion approximations have been studied in various models including those with customer abandonment and using a different scaling in state, namely, increasing the number of servers, which is called Halfin-Whitt regime. We will briefly discuss them.

%------------------------------------------------------ 
\subsection{Paper structure}
\label{sect:paper}
%------------------------------------------------------ 

This paper is made up by seven sections. In \sectn{single-class}, we introduce a stochastic model for the GJN, that is, generalized Jackson network, with a single-class of customers. This section have four subsections. 

In \sectn{primitive}, we start with a sample path level of modeling, then introduce stochastic assumptions for describing a GJN in \sectn{stochastic}. We are devoted to the single node case of the GJN, that is, the $GI/G/1$ queue, in Sections \sect{Kingman} and \sect{extension to queue length}. A main interest of these two subsections is to see how we should scale state and/or time to get a two-moments approximation for the stationary distribution. For this, we first review Kingman's \cite{King1962} approximation for the waiting time in \sectn{Kingman}, which is only scaled by state. In \sectn{extension to queue length}, we then consider the queue length process in continuous-time, and show that scaling is only required for state similarly to the waiting time case. They are turned out to be identical with the corresponding diffusion approximations. We conjecture that similar approximations can be obtained for the GJN with an arbitrary number of nodes.

In \sectn{diffusion}, we consider a diffusion approximation for a stochastic process. This section is divided into four subsections. In \sectn{SRBM}, we formally define the SRBM on a nonnegative multidimensional orthant, and derive a stationary equation to characterize the stationary distribution of SRBM, which is called a basic adjoint relationship, BAR for short (e.g., see \cite{HarrWill1987}). In \sectn{basic}, we present basic results for deriving a diffusion process as process limit. In \sectn{process limit}. those results will be used to derive the SRBM from the queue length processes of the GJNs under scaling state and time. We finally review the case that each node may have multiple servers in \sectn{extension to many servers}.

In \sectn{quality}, we consider quality of diffusion approximations for the stationary distribution of the GJN and some other characteristics. For the single node case, that is, a single queue, this approximation can be directly obtained for the stationary distribution as discussed in \sectn{extension to queue length}. In other words, a limiting process is not needed to get a weak limit of the stationary distribution. This also would be the case for the network case. However, in the literature, the SRBM is firstly obtained as a process limit, then the stationary distributions of the GJNs under diffusion scaling are shown to weakly converges to that of the SRBM. Surprisingly, this verification has been relatively recently done by Gamarnik and Zeevi \cite{GamaZeev2006} then by Budhiraja and Lee \cite{BudhLee2009} under weaker assumptions. We also review about diffusion approximations for moments and tail asymptotics of the stationary distribution.

In \sectn{multi-class}, we consider a multi-class queueing network. That is, types are attributed to customers, and nodes are exclusively attributed to types. Thus, the type of a customer uniquely specifies its service node. Each type at each node has its own primitive data, inter-arrival times (if it specifies exogenous arrivals) and service times. Customers may change their types when their services are completed. We assume that this change of types is Markovian. This network can not be handled by single-class networks in general because there may be more than one types of customers at a single node. In the nineties, multi-class queueing networks were highlighted because such models are found useful in production lines and unconventional stability conditions have been reported (e.g., see \cite{HarrNguy1995,Kuma1995}). These facts in addition to analytical difficulty stimulated people to study diffusion approximation for a multi-class queueing network. Harrison and Nguyen \cite{HarrNguy1990,HarrNguy1993} proposed a grand design for such diffusion approximation with some conjectures. Right after that, Dai and Wang \cite{DaiWang1993} negatively answered to those conjectures.

For a multi-class queueing network, a problem is under what conditions diffusion approximation exists in what sense. Harrison \cite{Harr1995} suggested to consider non-conventional diffusion approximation, but there seems no good alternative at least for open queueing networks. Thus, we here only consider the conventional diffusion approximation by diffusion scaling and SRBM. As for the conditions, state-space collapse and the completely-$\sr{S}$ condition for reflection matrix $R$ have been recognized as key features. We review diffusion approximation for a multi-class queueing network in this line.

In \sectn{Halfin-Whitt}, we review another type of diffusion approximations. Consider a $GI/G/s$ queue. For each fixed $s$, we can employ the standard diffusion approximation for this model. However, it ignores information on the number of busy servers, equivalently, the number of idle servers, so this approximation may be too rough for large $s$. This motivated Halfin and Whitt \cite{HalfWhit1981} to increase both of the arrival rate and the number $s$ of servers for a $GI/M/s$ queue while keeping the service time distribution, which is exponential in this model. This approximation also assumes heavy traffic condition, and called Halfin-Whitt regime. Its limiting process is also a diffusion process but state-dependent. This approximation and its variants have been largely studied because they can be used for call center application.

We finally give some concluding remarks in \sectn{concluding}.

%------------------------------------------------------ 
\section{Single-class queueing networks and heavy traffic approximation}
\label{sect:single-class}
%------------------------------------------------------ 

We first introduce a stochastic process for a single-class queueing model. Here, a single-class means that there is only one customer type, in other words, customers have no type. We assume that the model is open, and has exogenous arrivals. The system is composed of a finite number of nodes which have single waiting lines with possibly multiple servers, where waiting capacity is assumed to be infinite. Service times are attached at each servers. We assume that customers are served in the first-come first-served manner. If there are multiple nodes, customers may be routed to other nodes. 

\subsection{Primitive data and two characteristics}
\label{sect:primitive}

In this section, we assume that each node has a single server for simplicity. We describe this model starting at time $0$ by the following notations. Let $J \equiv \{1,2,\ldots,d\}$ be the set of nodes, where $d$ is a positive integer. For node $i \in J$, the sequence of service times $\{T_{si}(\ell); \ell \in \dd{N}\}$ is attached, where $\dd{N}$ is the set of natural numbers. Denote the set of nodes which have exogenous arrivals by  $J_{e}$, and, for $i \in J_{e}$, let $T_{ei}(\ell)$ be the inter-arrival time between the $(\ell-1)$-th and $\ell$-th exogenous arrivals at node $i$ for $\ell \in \dd{N}$. Note that $T_{ei}(\ell)$ may be $0$, which means a batch arrival. Let $\dd{Z}_{+}$ is the set of all nonnegative integers. For $n \in \dd{Z}_{+}$, let $\Phi_{ij}(n)$ be the number of customers who are routed to node $j$ among the first $n$ departures from node $i$, where $\Phi_{ij}(0) = 0$. This $j$ may be $0$, which means that customers leave the network. Thus,
\begin{align}
\label{eqn:primitive 1}
  \{T_{ei}(\ell); \ell \in \dd{N}\}, \quad i \in J_{e}, \quad \{T_{si}(\ell); \ell \in \dd{N}\}, \quad \{\Phi_{ij}(n); n \in \dd{Z}_{+},  j \in J \cup \{0\} \}, \quad i \in J,
\end{align}
specify this queueing network. We call them primitive data. To exclude trivial case, we assume that $J_{e} \ne \emptyset$.

For $i \in J_{e}$ and $i' \in J$, let
\begin{align*}
  V_{ei}(n) = \sum_{\ell=1}^{n} T_{ei}(\ell), \qquad V_{si'}(n) = \sum_{\ell=1}^{n} T_{si'}(\ell), \qquad n \in \dd{N},
\end{align*}
then the first two of the primitive data in \eqn{primitive 1} can be replaced by
\begin{align}
\label{eqn:primitive 2}
  \{V_{ei}(n); n \in \dd{N}\}, \quad i \in J_{e}, \qquad \{V_{si}(n); n \in \dd{N}\}, \quad i \in J.
\end{align}
We denote the counting processes for arrivals and service completions at node $i$ by
\begin{align}
\label{eqn:primitive 3}
  N_{ei}(t) \equiv \sup \{n; V_{ei}(n) \le t\}, \qquad N_{si}(t) \equiv \sup \{n; V_{si}(n) \le t\},
\end{align}
where $N_{ei}(t) \equiv 0$ for $i \in J \setminus J_{e}$. Then, \eqn{primitive 3} is equivalent to \eqn{primitive 2}, so can be used for the primitive data as well. The advantage of \eqn{primitive 3} is that the nodes which have no arrivals are not necessarily distinguished. Namely, we have just put $N_{ei}(t) \equiv 0$ for such nodes. It should be noted that $N_{ei}(t)$ counts exogenous arrivals at node $i$ but $N_{si}(t)$ may not count departures from node $i$. The latter does so if the server is always busy up to time $t$.

So far, we have assumed nothing about probability laws on those primitives, but we can define queue length and workload in sample path sense. Let $L_{i}(t)$ be the number of customers at node $i$ at time $t$, and $U_{i}(t)$ be the total busy time of the server at node $i$ up to time $t$. Then, it is easy to see that, for $t \ge 0$ and $i \in J$,
\begin{align}
\label{eqn:L 1}
 & L_{i}(t) = L_{i}(0) + N_{ei}(t) + \sum_{j \in J} \Phi_{ji}(N_{sj}(U_{j}(t))) - N_{si}(U_{i}(t)),\\
\label{eqn:U 1}
 & U_{i}(t) = \int_{0}^{t} 1(L_{i}(u) \ge 1) du,
\end{align}
because $N_{si}(U_{i}(t))$ is the number of departing customers from node $i$ in the time interval $[0,t]$, where $1(\cdot)$ is the indicator function of the statement ``$\cdot$'', which takes 1 if the statement is true, and vanishes otherwise. Conversely, it is not hard to see that $\{L_{i}(t); t \ge 0\}$ and $\{U_{i}(t); t \ge 0\}$ for $i \in J$ are uniquely determined by \eqn{L 1} and \eqn{U 1}. Hence, we can define them by those equations.

Since $N_{ei}(t) + \sum_{j \in J} \Phi_{ji}(N_{sj}(U_{j}(t)))$ is the total number of customers arrived at node $i$ by time $t$, the workload $W_{i}(t)$ at node $i$ can be defined as
\begin{align}
\label{eqn:WL 1}
 & W_{i}(t) = W_{i}(0) + V_{si}\Big(N_{ei}(t) + \sum_{j \in J} \Phi_{ji}(N_{sj}(U_{j}(t)))\Big) - U_{i}(t).
\end{align}
Let $\vc{L}(t)$ and $\vc{W}(t)$ be the vectors whose $i$-th entries are $L_{i}(t)$ and $W_{i}(t)$, respectively. Thus, we have obtained the mappings from the primitive data and the initial data $\vc{L}(0)$ and $\vc{W}(0)$ to $\vc{L}(\cdot) \equiv \{\vc{L}(t); t \ge 0\}$ and $\vc{W}(\cdot) \equiv \{\vc{W}(t); t \ge 0\}$. These mappings are called reflection mappings. In what follows we mainly consider the queue length process $\vc{L}(\cdot)$ because similar arguments can be applied for $\vc{W}(\cdot)$ with help of $\vc{L}(\cdot)$.

\subsection{Stochastic assumptions for GJN}
\label{sect:stochastic}

We now introduce the following stochastic assumptions on the primitive data \eqn{primitive 1}.
\begin{itemize}
\item [(\sect{single-class}a)] $\{T_{ei}(\ell); \ell \in \dd{N}\}$ with $i \in J_{e}$, $\{T_{si}(\ell); \ell \in \dd{N}\}$ and $\{\Phi_{ij}(n); n \in \dd{Z}_{+},  j \in J \cup \{0\} \}$ with $i \in J$ are independent.
\item [(\sect{single-class}b)] $\{T_{ei}(\ell); \ell \in \dd{N}\}$ for $i \in J_{e}$ and $\{T_{si}(\ell); \ell \in \dd{N}\}$ for $i \in J$ are sequences of $i.i.d.$ random variables taking values in $\dd{R}_{+}$, and, for each $i$, $T_{ei}(\ell)$ and $T_{si}(\ell)$ have finite and positive first moments, which are denoted by $m_{ei}$ and $m_{si}$, respectively. Let $\lambda_{ei} = m_{ei}^{-1}$ for $e \in J_{e}$ and $\lambda_{ei} = 0$ for $i \in J \setminus J_{e}$. Thus, $\lambda_{ei}$ is the exogenous arrival rate at node $i \in J$. Similarly, the service rate at node $i$ is denoted by $\lambda_{si} \equiv m_{si}^{-1}$.
\item [(\sect{single-class}c)] $\{\Phi_{ij}(n); n \in \dd{Z}_{+},  j \in J \cup \{0\} \}$ is a discrete-time $(d+1)$-dimensional counting process with $i.i.d.$ increments for each $i \in J$, and the increment $\Delta \Phi_{ij}(n) \equiv \Phi_{ij}(n) - \Phi_{ij}(n-1)$ takes values 0 or 1, and has the distribution given by
\begin{align*}
  \dd{P}(\Delta \Phi_{ij}(n) = 1) = p_{ij} \ge 0, \qquad j \in J.
\end{align*}
Let $P$ be the $d$-dimensional square matrix whose $(i,j)$ entry is $p_{ij}$ for $i,j \in J$, and let $\ol{P}$ be the $(d+1)$-dimensional square matrix whose $(i,j)$ entry is $p_{ij}$ for $i,j \in J \cup \{0\}$, where $p_{0j} = \lambda_{ej}(\sum_{i=1} \lambda_{ei})^{-1}$ for $j \in J$ and $p_{00} = 0$. We assume that $\ol{P}$ is stochastic and irreducible, which means that customers will leave the network in finite times with probability one.
\end{itemize}
The $d$-node queueing network satisfying those three sets of the conditions is called a generalized Jackson network, GJN for short. 

Under the assumptions (\sect{single-class}a)--(\sect{single-class}c), we can construct a Markov process for $\vc{L}(\cdot)$ supplemented by auxiliary characteristics. To this end, let $R_{ei}(t)$ and $R_{si}(t)$ be the residual times for the next exogenous arrival and service completion, respectively, at node $i$ at time $t$, where $R_{ei}(t) \equiv 0$ for $i \in J \setminus J_{e}$ and $R_{si}(t) = 0$ if the server at node $i$ is idle at time $t$. They are formally defined as
\begin{align*}
  R_{ei}(t) = V_{ei}(N_{ei}(t)+1) - t, \quad R_{si}(t) = (V_{si}(N_{si}(U_{i}(t))+1) - U_{i}(t))1(L_{i}(t) \ge 1).
\end{align*}
Denote the vectors whose $i$-th entries are $R_{ei}(t)$ and $R_{ei}(t)$ by $\vc{R}_{e}(t)$ and $\vc{R}_{s}(t)$. Then, it is not hard to see that $\{(\vc{L}(t),\vc{R}_{e}(t), \vc{R}_{si}(t)); t \ge 0\}$ is a Markov process. Furthermore, it is a strong Markov process (see, e.g., \cite{Dai1995}). We refer to this process as a supplemented Markov process for the GJN.

We next consider the stability of this Markov process, that is, the condition under which the stationary distribution exists. Recall that $\lambda_{ei}$ is the exogenous arrival rate at node $i \in J$. Consider the set of the following traffic equations for variables $\lambda_{ai}$ with $i \in J$. 
\begin{align}
\label{eqn:traffic 1}
  \lambda_{aj} = \lambda_{ej} + \sum_{i \in J} \lambda_{ai} p_{ij}, \qquad j \in J.
\end{align}
Denote the $d$-dimensional vectors whose entries are $\lambda_{ei}$ and $\lambda_{ai}$ by $\vc{\lambda}_{e}$ and $\vc{\lambda}_{a}$, respectively. Since $I-P$ is invertible by (\sect{single-class}c), the solution $\vc{\lambda}_{a}$ for the equations \eqn{traffic 1} is given by
\begin{align*}
  \vc{\lambda}_{a} = (I-P^{\rs{t}})^{-1} \vc{\lambda}_{e},
\end{align*}
where $P^{\rs{t}}$ is the transposition of $P$. Because $J_{e} \ne \emptyset$ and the irreducibility of $\ol{P}$ in (\sect{single-class}c), $\vc{\lambda}_{a}$ is positive. Let $\lambda_{si} = m_{si}^{-1}$, and denote the $d$-dimensional vector whose $i$-th entry is $\lambda_{si}$ by $\vc{\lambda}_{s}$. Then, it is well known that the supplemented Markov process is stable if and only if
\begin{align}
\label{eqn:stability 1}
  (I-P^{\rs{t}})^{-1} \vc{\lambda}_{e} < \vc{\lambda}_{s},
\end{align}
where the inequality of vectors represents entry-wise inequalities (see, e.g. \cite{Sigm1990}).

\subsection{Kingman's two-moment approximation in heavy traffic}
\label{sect:Kingman}

In this and next sections, we consider the $GI/G/1$ queue, as a pilot example for diffusion approximation. As we mentioned in \sectn{introduction}, Kingman \cite{King1962} is the first one to get the approximation which corresponds with diffusion approximation. He directly considered the stationary distribution of the waiting times in the $GI/G/1$ queue, and derived a limiting distribution under a certain scaling, using the first two moments of the interarrival and service times. In his formulation, no continuous-time parameter is involved, which is different from those in the previous two sections. Nevertheless, his approximation agrees with that obtained from diffusion approximation.

Here two questions arise. Is a similar approximation possible for the stationary distribution of the queue length process $L(t)$ ? If so, how it is related to diffusion approximation~? To address those questions for the $GI/G/1$ queue, which is the single node case of GJN, we first summarize Kingman's \cite{King1962} approximation for the stationary waiting time distribution, and make clear how approximation is produced from the limiting distribution. We then answer the questions in the next section.
 
We use the same notations of the GJN introduced in the previous two sections, but drop the index specifying a node for simplicity. For example, $L(t), R_{e}(t), R_{s}(t)$ and $W(t)$ stand for $L_{1}(t), R_{e1}(t), R_{s1}(t)$ and $W_{1}(t)$. We further drop the time index $t$ for those characteristics in the steady state like $L, R_{e}, R_{s}$ and $W$.

Let $D(n)$ be the waiting time, also called delay, of the $n$-th arriving customer. Obviously, we have
\begin{align*}
  D(n) = W((T_{e}(1) + T_{e}(2) + \ldots + T_{e}(n))-), \qquad n = 1,2,\ldots.
\end{align*}
However, this expression is not so useful because continuous-time process $W(t)$ is involved. Instead of it, Kingman directly considers the stationary distribution, using so called Lindley's equation. Assume that $\rho \equiv \lambda_{e}/\lambda_{s} < 1$, and let $D$ be a random variable subject to the stationary distribution of $D(n)$, then the distributional invariance of one step transition in the steady state yields Lindley's equation:
\begin{align}
\label{eqn:D 1}
  D \overset{w}{=} \max(0, D+T_{s}-T_{e}),
\end{align}
where $D, T_{s}, T_{e}$ are independent, and $\overset{w}{=}$ stands for the equality in distribution. In this formulation, no continuous-time process is needed.

Kingman's \cite{King1962} idea is to scale the waiting time $D$ in the steady state by
\begin{align*}
  \alpha = \frac {2(\lambda_{s} - \lambda_{e})} {\lambda_{e}^{2}(\sigma^{2}_{e} + \sigma^{2}_{s})},
\end{align*}
and to prove that $\alpha D$ weakly converges to the exponential distribution with unit mean, which is the stationary distribution of a one-dimensional SRBM, as $\alpha \downarrow 0$, under the finiteness conditions on some moments of $T_{e}$ and $T_{s}$. See Theorem 1 of K\"{o}llerstr\"{o}m \cite{Koll1974} for a comprehensive proof. Based on this convergence, Kingman \cite{King1961a} suggests the following approximation for the stationary waiting time $D$.
\begin{align}
\label{eqn:D approximation 1}
  \dd{P}(D \le x) \approx 1 - \exp\Big(- \frac {2(\lambda_{s} - \lambda_{e})} {\lambda_{e}^{2} (\sigma^{2}_{e} + \sigma^{2}_{s})} x \Big), \qquad x \ge 0.
\end{align}
Here, it is notable to distinguish the following two steps.
\begin{mylist}{3}
\item [(S1)] Prove a weak convergence, that is, convergence in distribution,  under appropriate changes of primitive data.
\item [(S2)] Derive an approximation from (S1).
\end{mylist}
(S1) is a purely mathematical issue, while (S2) is not more than proposal. Nevertheless, we have a clear rule for (S2) as we observed for the above $D$. Namely, the scaling factor is moved into the stationary distribution for approximation. When we are discussing about diffusion approximation, we often only consider (S1). This is because (S2) is automatically done according to the rule.

Kingman's as well as K\"{o}llerstr\"{o}m's arguments are very clear, but lack the idea of a sequence of models for the convergence. Such a sequence makes clear how the primitive data are changed. For this, we here introduce a sequence of the $GI/G/1$ queues indexed by $n$, where the primitive data of the $n$-th system are denoted by $\{T_{e}^{(n)}(\ell)\}$ and $\{T_{s}^{(n)}(\ell)\}$. Denote the stationary waiting time of the $n$-th queue by $D^{(n)}$. Define $\alpha^{(n)}$ as
\begin{align*}
  \alpha^{(n)} = \frac {2(\lambda^{(n)}_{s} - \lambda^{(n)}_{e})} {(\lambda^{(n)}_{s} \sigma^{(n)}_{s})^{2} + (\lambda^{(n)}_{e} \sigma^{(n)}_{e})^{2}}, \qquad n=1,2,\ldots.
\end{align*}
Then, Theorem 1 of K\"{o}llerstr\"{o}m \cite{Koll1974} reads as
\begin{align}
\label{eqn:W convergence 1}
  \lim_{n \to \infty} \dd{P}(\alpha^{(n)} D^{(n)} \le x) = 1 - e^{-x}, \qquad x \ge 0,
\end{align}
if the following conditions are satisfied.
\begin{align}
\label{eqn:1-node asym c1}
 & \lim_{n \to \infty} \alpha^{(n)} = 0, \qquad \alpha^{(n)} > 0, \quad n \ge 1,\\
\label{eqn:1-node asym c2}
 & \inf_{n \ge 1} ((\sigma^{(n)}_{s})^{2} + (\sigma^{(n)}_{e})^{2}) > 0,\\
\label{eqn:1-node asym c3}
 & \sup_{n \ge 1} \dd{E}((T^{(n)}_{e})^{2+\delta} + (T^{(n)}_{s})^{2+\delta})) < \infty, \qquad \mbox{for some } \delta > 0,
\end{align}

Because of \eqn{1-node asym c3}, $(\sigma^{(n)}_{s})^{2} + (\sigma^{(n)}_{e})^{2}$ are uniformly bounded in $n$, and therefore \eqn{1-node asym c1} is equivalent to
\begin{align}
\label{eqn:1-node asym c1a}
   \lim_{n \to \infty} (\lambda^{(n)}_{s} - \lambda^{(n)}_{e}) = 0,
\end{align}
which is further equivalent to $\rho^{(n)} \equiv \lambda^{(n)}_{e}/\lambda^{(n)}_{s} \to 1$ if $\lambda^{(n)}_{s}$ is uniformly away from $0$. Thus, we may divide (S1) into further two steps.
\begin{mylist}{3}
\item [(S1-1)] Set up a sequence of models indexed by positive integer $n$.
\item [(S1-2)] Scale the state space of the $n$-th system, and derive the limiting distribution.
\end{mylist}

In particular, if we choose $\lambda_{u}^{(n)}$ and $\sigma_{u}^{(n)}$ for $u=e,s$ in such a way that
\begin{align}
\label{eqn:1-node asym c4}
  \lim_{n \to \infty} \lambda^{(n)}_{u} = \lambda_{u} > 0, \qquad \lim_{n \to \infty} \sigma^{(n)}_{u} = \sigma_{u} > 0, \qquad u=e,s,
\end{align}
where $\lambda_{u}$ and $\sigma_{u}$ are positive constants, then \eqn{W convergence 1} becomes
\begin{align}
\label{eqn:W convergence 2}
  \lim_{n \to \infty} \dd{P}((1-\rho^{(n)}) D^{(n)} \le x) & = \lim_{n \to \infty} \dd{P}\Big(\alpha^{(n)} D^{(n)} \le \frac {\alpha^{(n)}}{1-\rho^{(n)}} x\Big) \nonumber\\
  & = 1 - \exp\Big(-\frac {2} {\lambda_{e}(\sigma_{e}^{2} + \sigma_{s}^{2})} x\Big), \qquad x \ge 0,
\end{align}
since $\lambda_{e} = \lambda_{s}$. This shows that, in (S1-2), $D^{(n)}$ can be scaled by $1-\rho^{(n)}$. 

\subsection{Extension of Kingman's approximation to the queue length}
\label{sect:extension to queue length}
 
We next consider diffusion approximation for the stationary distribution of $L(t)$, which is the queue length including customers being served, that is, the number of customers in system. In this case, continuous-time is involved, and one may think that it would be easier to consider scaling a process rather than scaling a stationary distribution. This may be the reason why diffusion approximation on the queue length has been much less studied through the stationary equation, while process limits have been well studied as we shall see in \sectn{process limit}.

We here challenge diffusion approximation for $L(t)$ through the stationary equation, and like to answer whether or not time scaling is really needed. Intuitively, no time scaling is needed because the stationary distribution is independent of time. Kingman's approximation supports this. However, this has not yet been confirmed, while process limit mandatorily requires state and time scaling to get a diffusion process.

To attack this problem, we use slightly stronger assumptions on the tails of the inter-arrival time $T_{e}$ and the service time $T_{s}$. We here recall the simplified notation system introduced in the previous section, which drops the index for a node. Let $\widetilde{F}_{e}$ and $\widetilde{F}_{s}$ be the moment generating functions of $T_{e}$ and $T_{s}$, respectively. That is,
\begin{align*}
  \widetilde{F}_{e}(v) = \dd{E}(e^{v T_{e}}), \qquad \widetilde{F}_{s}(v) = \dd{E}(e^{v T_{s}}), \qquad v \in \dd{R}.
\end{align*}
Assume the following light tail conditions.
\begin{itemize}
\item [(\sect{single-class}d)] There exist $\delta_{e}, \delta_{s} > 0$ such that
\begin{align}
\label{eqn:1-node condition 0}
 & \widetilde{F}_{u}(v) < \infty, \quad v < \delta_{u}, \qquad \lim_{v \uparrow \delta_{u}} \widetilde{F}_{u}(v) = \infty, \qquad u = e,s.
\end{align}
\end{itemize}
These conditions are a bit stranger than what we actually need, but typically assumed to study large deviations (e.g., see \cite{GlynWhit1994a}). 

For $\theta, \xi \le 0$, $\eta < \delta_{e}$ and some function $g(\theta)$ which will be determined later, let
\begin{align*}
  X(t) = \left\{\begin{array}{ll}
  g(\theta) e^{\eta R_{e}(t)}, \quad & L(t) = 0,\\
  e^{\theta L(t) + \eta R_{e}(t) + \xi R_{s}(t)}, \quad & L(t) \ge 1,
  \end{array} \right.
\end{align*}
and let $N_{d}$ be the counting process for departing customers. Since $X(t)$ has jumps at instants when $N_{e}(t)$ or $N_{d}(t)$ is increased, it is obvious to see that
\begin{align}
\label{eqn:RCL path}
  X(t) = X(0) & + \int_{0}^{t} X'(u)du + \int_{0}^{t} \Delta X(u) (N_{e}(du)+N_{d}(du)),
\end{align}
where $\Delta X(u) = X(u+) - X(u-)$. Let $\dd{E}_{e}, \dd{E}_{d}$ be the expectations under the Palm distributions with respect to $N_{e}, N_{d}$. Namely, for $u=e, d$ and a continuous function $f$ on $\dd{R}$,
\begin{align}
\label{eqn:Palm 1}
  \dd{E}_{u}(f(X(0))) = \lambda_{u}^{-1} \dd{E}\Big( \int_{0}^{1} f(X(u)) N_{u}(ds)\Big),
\end{align}
as long as the expectations are well defined, where $\lambda_{d} = \lambda_{e}$.
  
  Assume that $(L(0), R_{e}(0), R_{s}(0))$ is subject to the stationary distribution, which is also assumed to exist. Then, taking the expectation of both sides of \eqn{RCL path} yields
\begin{align*}
  0 = \dd{E}(X'(0)) + \lambda_{e} \dd{E}_{e}(\Delta X(0))  + \lambda_{d} \dd{E}_{d}(\Delta X(0)),
\end{align*}
  where the expectations are all finite for $\theta, \xi \le 0$, $\eta < \delta_{e}$. This formula is called a rate conservation law in the literature (e.g., see \cite{Miya1994}). Hence, letting
\begin{align*}
 & \varphi(\theta, \eta, \xi) = \dd{E}(e^{\theta L(0) + \eta R_{e}(0) + \xi R_{s}(0)}), \quad \varphi_{0}(\eta) = \dd{E}(e^{\eta R_{e}(0)}1(L(0)=0)), \\
 & \varphi_{e}(\theta, \xi) = \dd{E}_{e}(e^{\theta L(0-) + \xi R_{s}(0-)}), \qquad \varphi_{e0-}(0) = \dd{P}_{e}(L(0-)=0),\\
 & \varphi_{d}(\theta, \eta) = \dd{E}_{d}(e^{\theta L(0-) + \eta R_{e}(0-)}), \qquad \varphi_{d0+}(\eta) = \dd{E}_{d}(e^{\eta R_{e}(0)} 1(L(0+)=0)),
\end{align*}
we have the stationary equation for $\theta, \xi \le 0$, $\eta < \delta_{e}$:
\begin{align}
\label{eqn:1-node SE 1}
   -  (\eta + \xi) & \varphi(\theta, \eta, \xi) + (\eta (1-g(\theta)) + \xi) \varphi_{0}(\eta) \nonumber\\
   & + \lambda_{e} (e^{\theta} \widetilde{F}_{e}(\eta) - 1) \varphi_{e}(\theta, \xi) + \lambda_{e} (e^{\theta}  \widetilde{F}_{e}(\eta) ( \widetilde{F}_{s}(\xi)-1) + 1-g(\theta)) \varphi_{e0-}(0) \nonumber\\
  & + \lambda_{d} (e^{-\theta} \widetilde{F}_{s}(\xi) - 1) \varphi_{d}(\theta, \eta) + \lambda_{d} e^{-\theta} (g(\theta) - \widetilde{F}_{s}(\xi) ) \varphi_{d0+}(\eta) = 0.
\end{align}
  
Choose $\eta$, $\xi$ and $g(\theta)$ for $\theta \in \dd{R}$ in such a way that
\begin{align}
\label{eqn:1-node condition}
  e^{\theta} \widetilde{F}_{e}(\eta) = 1,\qquad e^{-\theta} \widetilde{F}_{s}(\xi) = 1, \qquad g(\theta) = e^{\theta}.
\end{align}
Then, they do exist by the assumption \eqn{1-node condition 0}, and \eqn{1-node SE 1} boils down to
\begin{align}
\label{eqn:1-node SE 2}
   (\eta + \xi) & \varphi(\theta, \eta, \xi) = (\eta (1-e^{\theta}) + \xi) \varphi_{0}(\eta), \qquad \theta, \xi \le 0, \; \eta < \delta_{e},
\end{align}
where $\varphi(\theta, \eta, \xi)$ and $\varphi_{0}(\eta)$ are finite because $\theta, \xi \le 0$ and
\begin{align*}
  \dd{E}(e^{\eta R_{e}}) = \lambda_{e} \dd{E}\Big( \int_{0}^{T_{e}} e^{\eta (T_{e} - x)} dx \Big) = \frac {\lambda_{e}}{\eta} \big(\widetilde{F}_{e}(\eta) - 1 \big) < \infty, \qquad \eta < \delta_{e},
\end{align*}
by the cycle formula, where the singularity of the right-hand term at the origin can be removed as it can be seen from the middle term.

We next note that $\eta$ and $\xi$ have nice representations as functions of $\theta$. For this, we define logarithmic moment generating functions for counting processes $N_{e}$ and $N_{s}$ as
\begin{eqnarray}
\label{eqn:1-node gamma}
  \gamma_{u}(\theta) = \lim_{t \to \infty} \frac 1t \log \dd{E}(e^{\theta N_{u}(t)}), \qquad u=e, s.
\end{eqnarray}
These functions play important roles in the theory of large deviations (e.g., see \cite{DembZeit1998}). We here note the following fact, which is immediate from Theorem 1 of \cite{GlynWhit1994a} (see also Lemma 4.1 of \cite{NeyNumm1987} in a more general context).

\begin{lemma}
\label{lem:1-node logarithmic}
  Under the assumption \eqn{1-node condition 0}, let $\eta(\theta)$ and $\xi(\theta)$ be the solutions $\eta$ and $\xi$ of \eqn{1-node condition}, then
\begin{align}
\label{eqn:1-node logarithmic}
  \eta(\theta) = - \gamma_{e}(\theta), \qquad \xi(\theta) = - \gamma_{s}(-\theta), \qquad \theta \in \dd{R}.
\end{align}
\end{lemma}

We here sketch a proof for the first equation of \eqn{1-node logarithmic} to see the reason why it is obtained (see \cite{GlynWhit1994a} for a complete proof). The second equation can be proved similarly. Define function $h$ for a given constant $\alpha$ as
\begin{align*}
  h(t) = \dd{E}\big( e^{\theta N_{e}(t) + \alpha t}), \qquad t \ge 0.
\end{align*}
It is easy to see that $h$ satisfies the following renewal equation.
\begin{align*}
  h(t) = e^{\alpha t} \dd{P}(T_{e} > t) + \int_{0}^{t} h(t-s) e^{\theta + \alpha s} d\dd{P}(T_{e} \le s).
\end{align*}
If we choose $\alpha$ such that $\widetilde{F}_{e}(\alpha) e^{\theta} = 1$, then this renewal equation has a proper interval distribution, and therefore we can see that $h(t)$ converges to a finite constant by the key renewal theorem as long as $\widetilde{F}_{e}(t)$ is finite around a neighborhood of $\alpha$. This is always the case under the assumption \eqn{1-node condition 0}. On the other hand, from the definition of $h$,
\begin{align*}
  \log h(t) = \log \dd{E}\big( e^{\theta N_{e}(t)}) + \alpha t.
\end{align*}
Dividing both sides of this equation by $t$ and letting $t$ to infinity, we have \eqn{1-node logarithmic}.

For scaling $L$, we expand $\gamma_{e}(\theta)$ and $\gamma_{s}(-\theta)$ as polynomials of $\theta$ around the origin.

\begin{lemma}
\label{lem:Taylor 1}
Under the same assumptions of \lem{1-node logarithmic}, we have, as $\theta \to 0$,
\begin{align}
\label{eqn:Taylor 1}
 & \gamma_{e}(\theta) = \lambda_{e} \theta + \frac 12 \lambda_{e}^{3} \sigma_{e}^{2} \theta^{2} + o(\theta^{2}), \\
\label{eqn:Taylor 2}
 & \gamma_{s}(-\theta) = - \lambda_{s} \theta + \frac 12 \lambda_{s}^{3} \sigma_{s}^{2} \theta^{2} + o(\theta^{2}).
\end{align}
\end{lemma}

This lemma is easily obtained from \eqn{1-node condition} and \eqn{1-node logarithmic}. For convenience of the reader, we prove it in \app{Taylor 1}.

We are now ready to consider a scaled limit for a sequence of the $GI/G/1$ queues, indexed by $n$. Assume that those queues are all stable, that is, $\rho^{(n)} < 1$ for all $n \ge 1$. As we did, characteristics of the $n$-th queue are indexed by superscript $^{(n)}$. By \lem{1-node logarithmic}, \eqn{1-node SE 2} for the $n$-th system can be written as
\begin{align*}
  -(\gamma^{(n)}_{e}(\theta) + \gamma^{(n)}_{s}(-\theta)) & \varphi^{(n)}(\theta, -\gamma^{(n)}_{e}(\theta), -\gamma^{(n)}_{s}(-\theta)) \nonumber\\
  & = (\gamma^{(n)}_{e}(\theta) (e^{\theta} - 1) - \gamma^{(n)}_{s}(-\theta)) \varphi^{(n)}_{0}(-\gamma^{(n)}_{e}(\theta)).
\end{align*}
Divide both sides of this equation by $\theta$ and applying \lem{Taylor 1}, we have
\begin{align}
\label{eqn:1-node SE 3}
 \Big( & (\lambda^{(n)}_{s} - \lambda^{(n)}_{e}) - \frac 12 ((\lambda^{(n)}_{e})^{3} (\sigma^{(n)}_{e})^{2} + (\lambda^{(n)}_{s})^{3} (\sigma^{(n)}_{s})^{2}) \theta + o(\theta) \Big) \nonumber\\
 & \hspace{25ex} \times \varphi^{(n)}(\theta, -\lambda^{(n)}_{e} \theta + o(\theta), \lambda^{(n)}_{s} \theta + o(\theta)) \nonumber\\
  & = \Big(\lambda^{(n)}_{e} (e^{\theta} - 1) + \lambda^{(n)}_{s} -  \frac 12 (\lambda^{(n)}_{s})^{3} (\sigma^{(n)}_{s})^{2} \theta + o(\theta) \Big)\nonumber\\
  & \hspace{25ex} \times (1 - \rho^{(n)}) \dd{E}(e^{(-\lambda^{(n)}_{e} \theta + o(\theta)) R^{(n)}_{e}(0)}| L(0)= 0),
\end{align}
where we have used the fact that $\dd{P}(L^{(n)} = 0) = 1 - \rho^{(n)}$.

Let us consider what is an appropriate scaling factor for $L^{(n)}$. Denote it by $\beta^{(n)} > 0$. That is, we scale $L^{(n)}$ as $\beta^{(n)} L^{(n)}$. This can be done by replacing $\theta$ by $\beta^{(n)} \theta$ in \eqn{1-node SE 3}. We first assume that $\beta^{(n)} \to 0$ as $n \to \infty$ because we hardly expect to have a reasonable limiting distribution otherwise. We then examine \eqn{1-node SE 3}, and easily see that we have a meaningful limiting distribution only when $\beta^{(n)}$ is proportional to $\lambda^{(n)}_{s} - \lambda^{(n)}_{e} > 0$. Under this limiting operation, it is preferable for the terms on the residual times to vanish. Namely, we like to verify the following conditions.
\begin{align}
\label{eqn:1-node extra 1}
 & \lim_{n \to \infty} \dd{E}(e^{- \beta^{(n)} \theta R^{(n)}_{e}}) = 1, \qquad \mbox{ for some } \theta < 0, \\
\label{eqn:1-node extra 2}
 & \lim_{n \to \infty} \beta^{(n)} R^{(n)}_{s} = 0 \; \mbox{ almost surely}.
\end{align}
Among these conditions, \eqn{1-node extra 2} is easily verified if the previous conditions \eqn{1-node asym c3} and \eqn{1-node asym c4} are satisfied because \eqn{1-node asym c3} implies that $\dd{E}(R^{(n)}_{s})$ is uniformly bounded in $n$. However, they may not be sufficient for \eqn{1-node extra 1}. So far, we here use the following simple assumption.
\begin{align}
\label{eqn:1-node extra 3}
  \mbox{The distribution of $T^{(n)}_{e}$ is independent of $n$.}
\end{align}

We now arrive at the following theorem, which corresponds to Kingman's approximation for the queue length in the steady state.

\begin{theorem}
\label{thr:1-node L asym}
  For a sequence of the stable $GI/G/1$ queues satisfying \eqn{1-node condition 0} for $\widetilde{F}_{s} = \widetilde{F}^{(n)}_{s}$, if the conditions \eqn{1-node asym c3}, \eqn{1-node asym c1a}, \eqn{1-node asym c4} and \eqn{1-node extra 3} hold and if $\rho^{(n)} < 1$ for $n \ge 1$, then
\begin{align}
\label{eqn:1-node L asym}
  \lim_{n \to \infty} \dd{P}( (1-\rho^{(n)})L^{(n)} \le x) = 1 - \exp\Big(-\frac {2} { \lambda_{e}^{2} (\sigma_{e}^{2} + \sigma_{s}^{2})} x\Big), \qquad x \ge 0.
\end{align}
\end{theorem}

We already explained major ideas to have this theorem, but its full proof requires careful arguments. Since they are technical, we defer them into \app{theorem L asym}. Note that the exponential distribution \eqn{1-node L asym} is identical with the stationary distribution of the SRBM, which will given in \cor{GJN 1}.

Thus, \thr{1-node L asym} concludes that time scaling is not needed to get the limiting distribution corresponding to that of SRBM. We here only considered a single node queue, but can expect a similar result for a general GJN although we have to work on the stationary equation instead of the stationary distribution because the latter is unknown for $d \ge 2$. Thus, we may have a completely different story about diffusion approximation for the stationary distribution of a queueing network in heavy traffic. As we will see in \sectn{quality}, this approach is closely related to a study on the tail decay rates of the stationary distributions.

In what follows, we will review another story, which requires scaling in both state and time.

\section{Diffusion approximation for a single class GJN}
\label{sect:diffusion}

We review diffusion approximation for a single class GJN introduced in Sections \sect{primitive} and \sect{stochastic} by an SRBM, that is, semi-martingale reflecting Brownian motion on a nonnegative orthant. This is now standard, and can be found in a text book (e.g., see \cite{ChenYao2001}). Nevertheless, it is good to see how the diffusion approximation is derived because it is a prototype for diffusion approximations for stochastic processes.

\subsection{From a single class GJN to SRBM}
\label{sect:SRBM}

We have observed that the supplemented Markov process may be useful to consider a limit of the stationary distributions. However, this is not the case for a limit of processes because an evolution in time needs to match scaling in state. Instead of it, the sample path characterizations \eqn{L 1} and \eqn{U 1} are known to be useful. To see this, let $R = I - P^{\rs{t}}$, and we rewrite \eqn{L 1} as
\begin{align}
\label{eqn:L 2}
 \vc{L}(t) = \vc{L}(0) + \vc{X}(t) + R \vc{Y}(t), \qquad t \ge 0,
\end{align}
where the $i$-th entry of $\vc{X}(t)$ is given by
\begin{align}
\label{eqn:X 1} 
  X_{i}(t) = & N_{ei}(t) - \lambda_{ei}t + \sum_{j \in J} \big(N_{sj}(U_{j}(t)) - \lambda_{sj} U_{j}(t) \big)  p_{ji} - (N_{si}(U_{i}(t)) - \lambda_{si} U_{i}(t)) \nonumber \\
  &  + \sum_{j \in J} \big(\Phi_{ji}(N_{sj}(U_{j}(t))) - N_{sj}(U_{j}(t)) p_{ji} \big) + \Big(\lambda_{ei} + \sum_{j \in J} \lambda_{sj} p_{ji} - \lambda_{si} \Big) t,
\end{align}
and the $i$-th entry of $\vc{Y}(t)$ is given by
\begin{align}
\label{eqn:Y i}
  Y_{i}(t) = \lambda_{si} (t - U_{i}(t)).
\end{align}
Clearly, $t - U_{i}(t)$ is the total idle time of the server at node $i$ by time $t$, and therefore $Y_{i}(t)$ is non-decreasing in $t$. Furthermore, $Y_{i}(t)$ increases only when $L_{i}(t) = 0$. These are precisely what are meant by \eqn{U 1}, and we have
\begin{align}
\label{eqn:Y 1}
  \int_{0}^{t} L_{i}(u) dY_{i}(u) = 0, \qquad i \in J.
\end{align}
Thus, \eqn{L 1} and \eqn{U 1} imply \eqn{L 2} and \eqn{Y 1}. Conversely, for a given $\vc{X}(\cdot)$, $(\vc{L}(\cdot), \vc{Y}(\cdot))$ is uniquely determined by \eqn{L 2} and \eqn{Y 1} as will be given in \thr{reflection map 1} below. Thus, we have the reflection mapping from $\vc{X}(\cdot)$ to $(\vc{L}(\cdot), \vc{Y}(\cdot))$, which is denoted by $(\phi(\vc{X}), \Psi(\vc{X}))$. This mapping will be used to study process limit for diffusion approximation.

A basic idea for diffusion approximation is to replace $\vc{X}(t)$ in \eqn{L 2} by a $d$-dimensional Brownian motion $\vc{B}(t)$ with a drift vector $\vc{\mu}$ satisfying $\vc{B}(0)=\vc{0}$. Namely, 
\begin{align}
\label{eqn:X 2}
  \vc{X}(t) = \vc{B}(t) + \vc{\mu} t.
\end{align}
Roughly speaking, $\{\vc{B}(t)\}$ is an independent incremental process with continuous sample path. The reader who is not familiar with Brownian motion may refer to Section 5.3 of \cite{ChenYao2001} for a short summary. Thorough introductions are found in \cite{Bill1999,Kall1997,StroVara1979}.

For this $\vc{X}(t)$, we rewrite \eqn{L 2} and \eqn{Y 1} for a nonnegative $d$-dimensional process $\vc{Z}(t)$ and a nondecreasing $d$-dimensional process $\vc{Y}(t)$ as
\begin{align}
\label{eqn:Z 1}
 & \vc{Z}(t) = \vc{Z}(0) + \vc{X}(t) + R \vc{Y}(t),\\
\label{eqn:Y 2}
 & \int_{0}^{t} Z_{i}(u) dY_{i}(u) = 0, \qquad i \in J, \qquad t \ge 0.
\end{align}

Let $D(\dd{R}^{d})$ be the set of functions from $[0,\infty)$ to $\dd{R}^{d}$ which are right continuous and have left-hand limit. In general, finding a solution $(\vc{Z}(\cdot), \vc{Y}(\cdot))$ in $D(\dd{R}^{2d})$ satisfying \eqn{Z 1} and \eqn{Y 2} for a given $\vc{Z}(0)$ and $\vc{X}(\cdot)$ with $\vc{X}(0) = \vc{0}$ is called a Skorohod problem. The following theorem answers this problem for $R = I - P^{\rs{t}}$.

\begin{theorem}[Theorem 1 of \cite{GamaZeev2006}, Theorem 7.2 of \cite{ChenYao2001}]
\label{thr:reflection map 1}
  Let $R = I - P^{\rs{t}}$. For each function $\vc{x}(\cdot) \in D(\dd{R}^{d})$, there exists a unique pair of nonnegative function $\vc{z}(\cdot)$ and a nondecreasing function $\vc{y}(\cdot)$ with $\vc{y}(0) = \vc{0}$ which satisfy
\begin{align}
\label{eqn:z 1}
 & \vc{z}(t) = \vc{x}(t) + R \vc{y}(t),\\
\label{eqn:y 1}
 & \int_{0}^{t} z_{i}(u) dy_{i}(u) = 0, \qquad i \in J, \qquad t \ge 0.
\end{align}
Furthermore, the mapping $\phi$ and $\Psi$ from $\vc{x}(\cdot)$ to $\vc{z}(\cdot) = \phi(\vc{x})$ and $\vc{y}(\cdot) = \Psi(\vc{x})$ are Lipschitz continuous. Namely, there exist constants $C_{1}$ and $C_{2}$ for each $t$ such that, for $H_{1} = \phi$ and $H_{2} = \Psi$,
\begin{align*}
 & \sup_{0 \le u \le t} \|H_{k}(\vc{x}_{1})(u) - H_{k}(\vc{x}_{2})(u) \| \le C_{k} \sup_{0 \le u \le t} \|\vc{x}_{1}(u) - \vc{x}_{2}(u))\|, \qquad k=1,2,
\end{align*}
for $\vc{x}_{1}(\cdot), \vc{x}_{2}(\cdot) \in D(\dd{R}^{d})$, where $\|\cdot\|$ is the Euclidean metric in the vector space $\dd{R}^{d}$.
\end{theorem}

If we choose $\vc{X}(\cdot)$ of \eqn{X 2} and take $\{\vc{Z}(0) + \vc{X}(t)\}$ for $\vc{x}(\cdot)$ in this theorem, we have a unique solution $(\vc{Z}(\cdot), \vc{Y}(\cdot))$ of \eqn{Z 1} and \eqn{Y 2}. Here, it is important for $R$ to be $I - P^{\rs{t}}$, which is always the case for GJN. However, some other models may require more general matrix $R$, which is indeed the case for a multi-class queueing network (see \sectn{multi-class}). Thus, we first define $(\vc{Z}(\cdot), \vc{Y}(\cdot))$ for a general square matrix $R$.

\begin{definition}
\label{den:SRBM}
Continuous $d$-dimensional processes $\{\vc{Z}(t); t \ge 0\}$ and $\{\vc{Y}(t); t \ge 0\}$ are called a semimartingale reflecting Brownian motion, SRBM for short, and a regulator, respectively, if there exists a family of probability measures $\{\dd{P}_{\vc{x}}; \vc{x} \in \dd{R}^{d}_{+}\}$ on some filtration $\{\sr{F}_{t}\}$ and the following conditions are satisfied for $\vc{Z}(0) = \vc{x}$ under $\dd{P}_{\vc{x}}$-almost surely for every $\vc{x} \in \dd{R}^{d}_{+}$.
\begin{itemize}
\item [(i)] $\{\vc{X}(t)\}$ is a $d$-dimensional Brownian motion with drift vector $\vc{\mu}$, that is, given by \eqn{X 2}, and $\{\vc{B}(t)\}$ is $\{\sr{F}_{t}\}$-martingale.
\item [(ii)] $\{\vc{Y}(t)\}$ is a $\{\sr{F}_{t}\}$-adapted nondecreasing process with $\vc{Y}(0) = \vc{0}$. 
\item [(iii)] $\{\vc{Z}(t)\}$ is a $\{\sr{F}_{t}\}$-adapted process.
\item [(iv)] \eqn{Z 1} and \eqn{Y 2} are satisfied.
\end{itemize}
\end{definition}

To use this definition, it is important to see when $\vc{Z}(\cdot)$ and $\vc{Y}(\cdot)$ uniquely exist. From \thr{reflection map 1}, we can see that, if $R = I - P^{\rs{t}}$, then they uniquely exist, which is independent of a choice of $\dd{P}_{\vc{x}}$. This existence problem was completely solved by Reiman and Williams \cite{ReimWill1988} for necessity and by Taylor and Williams \cite{TaylWill1993} for sufficiency. The answer is:
\begin{theorem}
\label{thr:SRBM 1}
  The SRBM $\vc{Z}(\cdot)$ uniquely exists if and only if $R$ is a completely-$\sr{S}$ matrix, which is a square matrix whose all principal submatrices are $\sr{S}$-matrices, where an $\ell$-dimensional square matrix $A$ is called an $\sr{S}$-matrix if there is an $\ell$-dimensional vector $\vc{v} > 0$ such that $A \vc{v} > \vc{0}$.
\end{theorem}

It is easy to see that $I - P^{\rs{t}}$ is a special case of a completely-$\sr{S}$ matrix because we can choose $(I - P^{\rs{t}})^{-1} \vc{1}$ for  the $\vc{v}$, where $\vc{1}$ is the $d$-dimensional vector whose entries are all units. 

In the remaining part of this section, we assume that $R$ is a completely $\sr{S}$-matrix, so the SRBM $\{\vc{Z}(t)\}$ exists. We denote the covariance matrix of the Brownian motion $\vc{B}(t)$ by $\Sigma \equiv \{\sigma_{ij}; i,j \in J\}$, and assume that $\Sigma$ is non-singular. Thus, $(\Sigma, \vc{\mu}, R)$ is the primitive data. Clearly, this SRBM is a Markov process with respect to filtration $\{\sr{F}_{t}\}$, and its dynamics is specified by It\^{o} integration formula. Namely, let $C^{2}_{b}(\dd{R}^{d})$ be the set of twice continuously differentiable and bounded functions from $[0,\infty)$ to $\dd{R}^{d}$, then, for $f \in C^{2}_{b}(\dd{R}^{d})$, computing the difference $f(\vc{Z}(u))-f(\vc{Z}(0))$ by the integration with respect to the integrator $d(\vc{X}(u) + R\vc{Y}(u))$ yields the following formula almost surely.
\begin{align}
\label{eqn:Ito 1}
f(\vc{Z}(u))-f(\vc{Z}(0)) =& \int_{0}^{t} (\nabla f(\vc{Z}(u)))^{\rs{t}} (\vc{\mu}\, du+d\vc{B}(u)) \nonumber\\
& + \int_{0}^{t} \sr{L}f(\vc{Z}(u))\,du + \int_{0}^{t} (\nabla f(\vc{Z}(u))^{\rs t}R \,d\vc{Y}(u), \qquad t \ge 0,
\end{align}
where 
\begin{align*}
 \nabla  f(\vc{x}) =\Big(\frac {\partial}{\partial x_{1}} f(\vc{x}),\ldots,\frac {\partial}{\partial x_{d}} f(\vc{x}) \Big)^{\rs t},\qquad
 \sr{L}f(\vc{x}) =\frac 12 \sum_{i,j \in J} \sigma_{ij} \frac {\partial^{2}}{\partial x_{i} \partial x_{j}} f(\vc{x}).
\end{align*}
If one is not familiar with It\^{o} integration, he may have a trouble with the second integration term in \eqn{Ito 1}. This term is so called an integration by quadratic variations, which arise from unbounded variation of the sample path of a Brownian motion. On the other hand, the integration with respect to $d\vc{B}(u)$ in the first integration is just a Riemann integral, whose expectation vanishes. We now have nice text books for It\^{o} integration. For example, Chung and Williams \cite{ChunWill1990} is recommended for first reading.

We next consider a stable SRBM, that is, an SRBM which has a stationary distribution. The stability condition has not yet been completely obtained in terms of the primitive data for $d \ge 4$ (see \cite{BramDaiHarr2010}). However, if $R$ is an $\sr{M}$-matrix, which is a matrix of the form $I - G$ for some nonnegative matrix $G$ such that $I - G$ is invertible, then the SRBM is stable if and only if
\begin{align}
\label{eqn:stability 2}
  R^{-1} \vc{\mu} < \vc{0}.
\end{align}
In particular, this is the stability condition for $R = I - P^{\rs{t}}$ since this $R$ is an $\sr{M}$-matrix.

We return to a general $R$ which is completely $\sr{S}$, and assume that the SRBM is stable. We denote its stationary distribution by $\pi$, and define measure $\nu_{i}$ on $\dd{R}_{+}^{d}$ as
\begin{align}
\label{eqn:nu i}
  \nu_{i}(A) = \dd{E}_{\pi}\Big( \int_{0}^{1} 1(\vc{Z}(u) \in A) dY_{i}(u)\Big), \qquad A \in \sr{B}(\dd{R}_{+}^{d}), \; i \in J,
\end{align}
where $\sr{B}(\dd{R}_{+}^{d})$ is the Borel $\sigma$-field on $\dd{R}_{+}^{d}$, and $\dd{E}_{\pi}$ stands for the expectation under $\pi$, that is, under the condition that $\vc{Z}(0)$ is subject to $\pi$. Obviously, $\nu_{i}$ is a finite measure because $\dd{E}(X_{i}(t) - X_{i}(0))$ is bounded for $t \in [0,1]$. In the literature (e.g., see \cite{Kall1997}), the measure $\nu_{i}$ is called a Palm measure with respect to the random measure generated by $Y_{i}(\cdot)$. Since $Y_{i}(t)$ does not increases for $Z_{i}(t) > 0$, this measure concentrates on the face $F_{i}$ given by
\begin{align*}
  F_{i} = \{\vc{x} \in \dd{R}_{+}^{d}; x_{i} = 0\}.
\end{align*}
Let us take the expectation of \eqn{Ito 1} under $\pi$, then it is immediate to have the part (a) of the following theorem.

\begin{theorem}
\label{thr:BAR}
Let $R$ be a completely $\sr{S}$-matrix. (a) If the SRBM has the stationary distribution $\pi$, then, for $f \in C_{b}^{2}(\dd{R}^{d})$,  
\begin{align}
\label{eqn:BAR 1}
  \int_{\dd{R}_{+}^{d}} \sr{L} f(\vc{x}) \pi(d\vc{x}) + \int_{\dd{R}_{+}^{d}} \br{\nabla f(\vc{x}), \vc{\mu}} \pi(d\vc{x}) + \sum_{i \in J} \int_{\dd{R}_{+}^{d}} \br{\nabla f(\vc{x}), R^{[i]}} \nu_{i}(d\vc{x}) = 0,
\end{align}
where $\br{\vc{a}, \vc{b}}$ is the inner product $\vc{a}^{\rs{t}} \vc{b}$ of vectors $\vc{a}, \vc{b} \in \dd{R}^{d}$, and $R^{[i]}$ is the $i$-th column of $R$. (b) Conversely, if there are a probability distribution $\pi$ on $(\dd{R}_{+}^{d}, \sr{B}(\dd{R}_{+}^{d}))$ and finite measures $\nu_{i}$ on $(F_{i}, \sr{B}(F_{i}))$ for $i \in J$, then $\pi$ is the stationary distribution of the SRBM, and \eqn{nu i} holds.
\end{theorem}
\begin{remark}
\label{rem:BAR}
  \eqn{BAR 1} for all $f \in C_{b}^{2}(\dd{R}^{d})$ is equivalent to
\begin{align}
\label{eqn:BAR 2}
   \Big(\frac 12 \br{\vc{\theta}, \Sigma \vc{\theta}} + \br{\vc{\mu}, \vc{\theta}}\Big) \tilde{\varphi}(\vc{\theta}) + \sum_{i \in J} \br{R^{[i]}, \vc{\theta}} \tilde{\varphi}_{i}(\vc{\theta}) = 0, \qquad \vc{\theta} \le \vc{0},
\end{align}
where $\tilde{\varphi}(\vc{\theta})$ and $\tilde{\varphi}_{i}(\vc{\theta})$ are the moment generating functions of the stationary distribution $\pi$ and the Palm measure $\nu_{i}$, respectively. \eqn{BAR 2} is immediate from \eqn{BAR 1} with $f(\vc{x}) = e^{\br{\vc{\theta}, \vc{x}}}$ for $\vc{x} \in \dd{R}_{+}^{d}$. The converse is also easy to see at least intuitively because Laplace transform uniquely determines a distribution and measure on the nonnegative orthant. A detailed proof for the converse is given in Appendix D of the arXiv version \cite{DaiMiyaWu2014b}.
\end{remark}

This equation \eqn{BAR 1} is called a basic adjoint relationship, BAR for short (e.g., see \cite{HarrWill1987}). Since \eqn{BAR 1} is nothing but a stationary equation, (b) must be true intuitively. However, its proof is technically quite involving (see \cite{DaiKurt1994,KangRama2014}).

\subsection{Basic tips for process limit}
\label{sect:basic}

For diffusion approximation, we will consider a sequence of stochastic processes which converges to some stochastic process. Recall that this limit is called a process limit. To formally define this convergence, we need topology or metric on the set of functions. We summarize some basic tips for them. The reader may skip this subsection if it is familiar with them. The materials here can be found in standard text books \cite{Bill1999,ChenYao2001,Kall1997,Whit2002}. The book chapter \cite{Glyn1990} is a good reference for OR readers.

Let $S$ be a complete and separable normed space with norm $\|\cdot\|_{S}$, and let $C(E, S)$ ($D(E, S)$) be the set of all functions from an interval $E \subset \dd{R}_{+}$ to $S$ which are continuous (are right-continuous and have left-hand limits, respectively). For $E=\dd{R}_{+}$, we omit $E$, so $C(\dd{R}_{+}, S)$ and $D(\dd{R}_{+}, S)$ are simply denoted by $C(S)$ or $D(S)$, which we already used for $C(\dd{R}^{d})$ or $D(\dd{R}^{d})$ in \sectn{SRBM}.

Stochastic processes which we consider are functions from the sample space $\Omega$ to either $C(S)$ or $D(S)$. In particular, a sample path of a Brownian motion belongs to $C(S)$ with $S = \dd{R}^{d}$, while a sample path for the GJN belongs to $D(S)$. Thus, we need topology or metric on both $C(S)$ and $D(S)$ to consider diffusion approximation. They can be defined through metrics on $C([0,t], S)$ and $D([0,t], S)$ for each $t > 0$. Let 
\begin{align*}
  \|x-y\|_{t} = \sup_{u \in [0,t]} \|x(u) - y(u)\|_{S}, \qquad x, y \in C([0,t],S) \mbox{ or } D([0,t],S),
\end{align*}
which is called a uniform norm. $C([0,t], S)$ typically uses this uniform norm as its metric. On the other hand, $D([0,t], S)$ uses $J_{1}$-metric, which is defined as
\begin{align*}
  d_{J_{1}}(x,y,t) = \inf_{\lambda \in \Lambda(t)} \max(\|x \circ \lambda - y\|_{t}, \|\lambda - \lambda_{0}\|_{t}), \qquad x, y \in D([0,t],S),
\end{align*}
where $\Lambda(t)$ is the set of strictly increasing and continuous function from $[0,t]$ to $[0,t]$, $x \circ \lambda(u) = x(\lambda(u))$ and $\lambda_{0}(u) = u$ for $u \in [0,t]$. As one can see from its definition, $J_{1}$-metric is to minimize discontinuity of functions by a suitable time scaling.

Let $X_{n}(\cdot) \equiv \{X_{n}(t); t \ge 0\}$ for $n=1,2,\ldots$ and $X(\cdot) \equiv \{X(t); t \ge 0\}$ be stochastic processes whose sample paths are in $D(S)$. Then, $X_{n}(\cdot)$ is said to converge to $X(\cdot)$ in distribution, which is denoted by $X_{n}(\cdot) \overset{w}{\longrightarrow} X(\cdot)$, if, for any bounded continuous function $f$ from $D(S)$ to $\dd{R}$,
\begin{align}
\label{eqn:convergence 1}
  \lim_{n \to \infty} \dd{E}(f(X_{n}(\cdot))) = \dd{E}(f(X(\cdot))).
\end{align}
We also refer to this convergence as weak convergence.

The following theorem is an immediate consequence of these definitions. 
\begin{theorem}[Continuous mapping theorem]
\label{thr:continuous}
Let $g$ be a continuous mapping from $D(S)$ to $D(S)$, and assume that $X_{n}(\cdot) \overset{w}{\longrightarrow} X(\cdot)$ for the sequence of stochastic processes $X_{n}(\cdot)$ and that of $X(\cdot)$ whose sample paths in $D(S)$, then we have
\begin{align*}
  g(X_{n})(\cdot) \overset{w}{\longrightarrow} g(X)(\cdot).
\end{align*}
\end{theorem}
Another useful tool is to replace convergence in distribution by almost surely convergence. For this, sample paths for the sequence of stochastic processes must be suitably chosen, but it is quite useful to prove the weak convergence of functionals of those stochastic processes. A random variable version of this theorem is not hard to prove (e.g., see \cite{Bill2000}), but its proof is quite technical for stochastic processes.

\begin{theorem}[Skorohod representation theorem, Theorem 6.7 of \cite{Bill1999}]
\label{thr:Skorohod}
  For the sequence of stochastic processes $X_{n}(\cdot)$ and that of $X(\cdot)$ whose sample paths are in $D(S)$, if $X_{n}(\cdot) \overset{w}{\longrightarrow} X(\cdot)$, then there exist stochastic processes $\tilde{X}_{n}(\cdot)$ and $\tilde{X}(\cdot)$ such that they have the same distributions as those of $X_{n}(\cdot)$ and $X(\cdot)$, respectively, and $\tilde{X}_{n}(\cdot)$ almost surely converges to $\tilde{X}(\cdot)$ as $n \to \infty$.
\end{theorem}

In general, verifying the convergence condition \eqn{convergence 1} is not easy. A typical situation is to use a finite dimensional convergence of $X_{n}(\cdot)$ in distribution. Namely, for each $k \ge 1$, each sequence $t_{1}, t_{2}, \ldots, t_{k} \in \dd{R}_{+}$ and each bounded continuous function $f_{k}$ from $S^{k}$ to $\dd{R}$,
\begin{align*}
  \lim_{n \to \infty} \dd{E}(f_{k}(X_{n}(t_{1}), X_{n}(t_{2}), \ldots, X_{n}(t_{k}))) = \dd{E}(f_{k}(X(t_{1}), X(t_{2}), \ldots, X(t_{k}))).
\end{align*}
Clearly, $X_{n}(\cdot) \overset{w}{\longrightarrow} X(\cdot)$ implies this finite dimensional convergence in distribution, but its converse is generally not true although the limiting distribution of $X_{n}(\cdot)$ is uniquely identified if it exists (e.g., see \cite{Glyn1990}). For the converse to be true, an extra condition is required. A typical condition for this is tightness. The sequence of $X_{n}(\cdot)$ is said to be tight if, for any $\epsilon > 0$, there exists a compact set $K$ of $D(S)$ such that
\begin{align*}
  \liminf_{n \to \infty} \dd{P}( X_{n}(\cdot) \in K) > 1- \epsilon.
\end{align*}
This tightness also may not be easy to verify. So, various sufficient conditions have been studied in the literature (e.g., see \cite{Bill1999,Kall1997}). Those are quite technical, and we will not get into their details. 

Another useful technique is a random change of time. 
\begin{lemma}[Lemma on page 151 of \cite{Bill1999}]
\label{lem:random change}
  Let $\{X^{(n)}(\cdot); n=1,2,\ldots\}$ be the sequence of processes whose sample paths are in $D(S)$, and let $A^{(n)}(\cdot)$ be the sequence of non-decreasing non-negative valued processes. If there are $X^{(n)}(\cdot) \in D(S)$ and nondecreasing $A^{(n)}(\cdot)$ such that $(X^{(n)}(\cdot), A^{(n)}(\cdot))  \overset{w}{\longrightarrow} (X(\cdot), A(\cdot))$ and if $\dd{P}(X(\cdot) \in C(S)) = 1$, then $X^{(n)}(\cdot) \circ A^{(n)}(\cdot) \overset{w}{\longrightarrow} X(\cdot) \circ A(\cdot)$, where the composition $f \circ g$ is defined as $f \circ g (x) = f(g(x))$.
\end{lemma}

We finally refer to process versions of the central limit theorems (called Lindeberg-Feller Theorems, e.g., see Theorem 7.2.1 of \cite{Chun2001}). They are all sources of diffusion approximation.
 
\begin{theorem}[An extended version of Donsker's theorem (e.g., see Problem 8.4 of \cite{Bill1999})]
\label{thr:Donsker}
  For each $n \ge 1$, let $\tau^{(n)}_{i}$, $\tau^{(n)}_{2}$, \ldots be independent and identically distributed random variables with finite mean $m^{(n)}$ and finite variance $(\sigma^{(n)})^{2}$. Assume that
\begin{align}
\label{eqn:sigma convergence}
 & \lim_{n \to \infty} \sigma^{(n)} = \sigma > 0,\\
\label{eqn:Lindeberg 3}
 & \sup_{n \ge 1} \dd{E}((\tau^{(n)})^{2+\delta}) < \infty, \qquad \exists \delta > 0,
\end{align}
and define $\widehat{X}^{(n)}(t)$ as
\begin{align*}
  \widehat{X}^{(n)}(t) = \frac 1{\sqrt{n} \sigma^{(n)}} \sum_{\ell=1}^{[nt]} (\tau^{(n)}_{\ell} - m^{(n)}),
\end{align*}
then $\widehat{X}^{(n)}(\cdot) \overset{w}{\longrightarrow} B(\cdot)$, where $[a]$ is the largest integer not greater than real number $a$, and $B(\cdot)$ is the standard Brownian motion.
\end{theorem}
\begin{remark}
\label{rem:Donsker}
  The conditions \eqn{sigma convergence} and \eqn{Lindeberg 3} guarantees for Lindeberg condition to hold (see Theorem 27.3 of \cite{Bill2000}).
\end{remark}

A similar scaling limit can be obtained for the counting process $N^{(n)}$, which is defined as
\begin{align*}
  N^{(n)}(t) = \sup \Big\{\ell \ge 0; \sum_{k=1}^{\ell} \tau^{(n)}_{k} \le t \Big\}.
\end{align*}
The next theorem is obtained from \thr{Donsker} using \lem{random change} (e.g., see \cite{Glyn1990} for details). 

\begin{theorem}[Counting process version of Donsker's theorem, Theorem 14.6 of \cite{Bill1999}]
\label{thr:counting}
 Assume that there exist positive and finite $m$ and $\sigma$ such that $\widehat{X}^{(n)}(\cdot) \overset{w}{\longrightarrow} B(\cdot)$ holds under the assumptions of \thr{Donsker}. Let $\lambda^{(n)} = 1/m^{(n)}$. Assume that
\begin{align*}
  \lim_{n \to \infty} \lambda^{(n)} = \lambda,
\end{align*}
and define
\begin{align*}
  \widehat{Z}^{(n)}(t) = \frac 1{\sqrt{n} \sigma^{(n)} (\lambda^{(n)})^{\frac 32}} \big( N^{(n)}(nt) - nt \lambda^{(n)} \big), \qquad t \ge 0,
\end{align*}
then $\widehat{Z}^{(n)}(\cdot) \overset{w}{\longrightarrow} B(\cdot)$.
\end{theorem}

From \lem{Taylor 1}, we can see that the variance $V(N^{(n)}(t)) \sim (\lambda^{(n)})^{3} (\sigma^{(n)})^{2} t$ as $t \to \infty$. This is the reason why $\sigma^{(n)} (\lambda^{(n)})^{\frac 32}$ appears in the normalizing factor above. For convenience of the reader, we directly verify this fact in \app{Donsker counting}.
 
\subsection{Process limit for diffusion approximation}
\label{sect:process limit}

We now consider a process limit of the joint queue length process $\{\vc{L}(t)\}$ of the GJN under suitable scaling in time and state. A key step for obtaining this limit is to change $X_{i}(t)$ of \eqn{X 1} to the following $\ul{X}_{i}(t)$ by replacing accumulated busy times $U_{i}(t)$ and $U_{j}(t)$ by deterministic time $t$.
\begin{align}
\label{eqn:X 3}
  \ul{X}_{i}(t) = & N_{ei}(t) - \lambda_{ei}t + \sum_{j \in J} \big(N_{sj}(t) - \lambda_{sj} t\big) (p_{ji} - \delta_{ij}) \nonumber \\
  & + \sum_{j \in J} \big(\Phi_{ji}(N_{sj}(t)) - N_{sj}(t) p_{ji} \big) + \Big(\lambda_{ei} + \sum_{j \in J} \lambda_{sj} p_{ji} - \lambda_{si} \Big) t,
\end{align}
where $\delta_{ij} = 1(i=j)$. We denote the vector whose $i$-th coordinate is $\ul{X}_{i}(t)$ by $\ul{\vc{X}}(t)$. 

One can guess that the process $\{\ul{\vc{X}}(t)\}$ converges to a Brownian motion under normalization. For this, we require additional conditions on the finiteness of variances. Denote generic random variables which have the same distributions as $T_{ei}(\ell)$ and $T_{si}(\ell)$ respectively by $T_{ei}$ and $T_{si}$. In what follow, we assume:
\begin{itemize}
\item [(\sect{diffusion}a)] $T_{ei}$ and $T_{si}$ have finite variances, which are denoted respectively by $\sigma_{ei}^{2}$ and $\sigma_{si}^{2}$, where we let $\sigma_{ei} = 0$ for $i \in J \setminus J_{e}$, and
\begin{align}
\label{eqn:non-trivial 1}
  \sum_{i \in J} (\sigma_{ei}^{2} + \sigma_{si}^{2}) > 0.
\end{align}
\end{itemize}
The condition \eqn{non-trivial 1} is just to exclude a trivial case in a simple way, but can be weakened as we will see later.

We next introduce the sequence of the GJN's for diffusion approximation. We index those GJN's by positive integer $n$. For the $n$-th GJN, all random elements as well as constant parameters are indexed by $n$ such as $T^{(n)}_{ei}(\ell)$, $T^{(n)}_{si}(\ell)$, $\lambda^{(n)}_{ei}$, $\vc{L}^{(n)}(t)$, $\vc{X}^{(n)}(t)$ and so on, while keeping the non-indexed notations for the original GJN satisfying the assumptions (\sect{single-class}a)--(\sect{single-class}c) and (\sect{diffusion}a). In particular, we will not change the routing matrix $P$, so $R$ is unchanged. We define the $n$-th GJN by
\begin{align}
\label{eqn:Tn 1}
 & T^{(n)}_{ei}(\ell) = T_{ei}(\ell), \qquad i \in J_{e}\\
\label{eqn:Sn 1}
 & T^{(n)}_{si}(\ell) = T_{si}(\ell) - m_{si} + \Big(1 - \frac 1 {\sqrt{n}} \Big) \lambda_{ai}^{-1}, \qquad i \in J,
\end{align}
keeping the other primitive data, where we recall that $\lambda_{ai} (\equiv \lambda^{(n)}_{ai})$ is the total arrival rate at node $i$, which is obtained by \eqn{traffic 1}. Then, it is easy to see that
\begin{align}
\label{eqn:SR 1}
  \sqrt{n} \big(\lambda^{(n)}_{ai} - \lambda^{(n)}_{si}\big) = - \lambda_{si} < 0, \qquad i \in J.
\end{align}
On the other hand, the variances of $T^{(n)}_{ei}$ and $T^{(n)}_{si}$ are unchanged, that is,
\begin{align*}
  \sigma^{(n)}_{ei} = \sigma_{ei}, \qquad \sigma^{(n)}_{si} = \sigma_{si}, \quad i \in J.
\end{align*}
Instead of those assumptions, one can choose $\lambda^{(n)}_{ei}, \lambda^{(n)}_{si}, \sigma^{(n)}_{ei}, \sigma^{(n)}_{si}$ as long as the following two sets of conditions are satisfied (e.g., see \cite{Reim1984}).
\begin{itemize}
\item [(\sect{diffusion}b)] $\lim_{n \to \infty} \lambda^{(n)}_{ei} = \lambda_{ei}$, $\lim_{n \to \infty} \lambda^{(n)}_{si} = \lambda_{si}$, $\lim_{n \to \infty} \sigma^{(n)}_{ei} = \sigma_{ei}$, $\lim_{n \to \infty} \sigma^{(n)}_{si} = \sigma_{si}$, and\\
$\lim_{n \to \infty} \sqrt{n} \big(\lambda^{(n)}_{ai} - \lambda^{(n)}_{si}\big) = c_{i}$ for some $c_{i} \in \dd{R}$.
\item [(\sect{diffusion}c)] $\sup_{n \ge 1} \dd{E}((T^{(n)}_{ei})^{2+\delta} 1(i \in J_{e}) + (T^{(n)}_{si})^{2+\delta}) < \infty$ for some $\delta > 0$.
\end{itemize}

In what follows we simply assume \eqn{Tn 1} and \eqn{Sn 1}. Obviously, (\sect{diffusion}b) is satisfied for this case, and (\sect{diffusion}c) is not needed because Lindeberg condition for the central limit theorem is automatically satisfied in this case. For the $n$-th GJN, we define diffusion scaling processes by
\begin{align*}
 & \widehat{\vc{L}}^{(n)}(t) = \frac 1{\sqrt{n}} \vc{L}^{(n)}(nt), \qquad \widehat{\vc{Y}}^{(n)}(t) = \frac 1{\sqrt{n}} \vc{Y}^{(n)}(nt), \\
 & \widehat{\vc{X}}^{(n)}(t) = \frac 1{\sqrt{n}} \vc{X}^{(n)}(nt), \qquad \widehat{\ul{\vc{X}}}^{(n)}(t) = \frac 1{\sqrt{n}} \ul{\vc{X}}^{(n)}(nt).
\end{align*}
We would like to show that $\{\widehat{\vc{X}}^{(n)}(t)\}$ has the same limiting process as that of $\{\widehat{\ul{\vc{X}}}^{(n)}(t)\}$. For this, we use the following fact.

\begin{lemma}[Proposition 4 of \cite{Reim1984}]
\label{lem:U 2}
$\ol{U}^{(n)}_{i}(t) \equiv n^{-1} U^{(n)}_{i}(nt)$ converges to $t$ in probability as $n \to \infty$.
\end{lemma}

By this lemma and \lem{random change}, we have
\begin{align}
\label{eqn:tilde X 1}
  \widehat{X}^{(n)}_{i}(t) = \frac 1{\sqrt{n}} \ul{X}^{(n)}_{i}(U^{(n)}_{i}(nt)) = \frac 1{\sqrt{n}} \ul{X}^{(n)}_{i}(n \ol{U}^{(n)}_{i}(t)) \simeq \ul{\widehat{X}}^{(n)}_{i}(t), \quad n \to \infty,
\end{align}
where $\simeq$ stands to have the same distribution asymptotically as $n \to \infty$.

Thus, for $\widehat{\vc{L}}^{(n)}(t)$ to weakly converges to an SRBM, it remains to show that
\begin{align*}
  \widehat{\ul{X}}^{(n)}_{i}(t) = & \frac 1{\sqrt{n}} \big(N^{(n)}_{ei}(nt) - \lambda^{(n)}_{ei} nt \big) + \sum_{j \in J} \frac 1{\sqrt{n}} \big(N^{(n)}_{sj}(nt) - \lambda^{(n)}_{sj} nt \big) (p_{ji} - \delta_{ji}) \nonumber \\
  & + \sum_{j \in J} \frac 1{\sqrt{n}} \big(\Phi_{ji}(N^{(n)}_{sj}(nt)) - N_{sj}^{(n)}(nt) p_{ji} \big) + \sqrt{n} \Big(\lambda^{(n)}_{ei} + \sum_{j \in J} \lambda^{(n)}_{sj} p_{ji} - \lambda^{(n)}_{si} \Big) t
\end{align*}
weakly converges to a $d$-dimensional Brownian motion with a constant drift term, where $\lambda^{(n)}_{ei} = \lambda_{ei}$ by \eqn{Tn 1}, but we keep the notation $\lambda^{(n)}_{ei}$ for arguments to be parallel for service times. Letting
\begin{align*}
  \widehat{N}^{(n)}_{ei}(t) \equiv \frac 1{\sqrt{n}} \big(N^{(n)}_{ei}(nt) - \lambda^{(n)}_{ei}nt\big),
\end{align*}
it follows from \thr{counting} that
\begin{align}
\label{eqn:ei to B1}
 \widehat{N}^{(n)}_{ei}(\cdot) \overset{w}{\longrightarrow} \sqrt{\lambda^{3}_{ei} \sigma_{ei}^{2}} B^{*}_{ei}(\cdot),
\end{align}
where $\{B^{*}_{ei}(t)\}$ is the standard Brownian motion. We here choose $\{B^{*}_{ei}(t)\}$ to be independent for $i \in J_{e}$.

Let $\widehat{\vc{N}}_{e}^{(n)}(t)$ is the random vector whose $i$-th coordinate is $\widehat{N}^{(n)}_{ei}(t)$. We similarly define $\widehat{\vc{N}}_{s}^{(n)}(t)$. It is immediate from \eqn{ei to B1} that
\begin{align}
\label{eqn:tilde N e}
  \widehat{\vc{N}}_{e}^{(n)}(\cdot) \overset{w}{\longrightarrow} \,\diag(\{\lambda^{\frac 32}_{ei} \sigma_{ei}\}) \vc{B}_{e}^{*}(\cdot),
\end{align}
where $\diag(\{a_{i}\})$ is the diagonal matrix whose $i$-th diagonal entry is $a_{i}$, and $\vc{B}_{e}^{*}(\cdot)$ is the $d$-dimensional standard Brownian motion, which means that its covariance matrix is the identity matrix. Note that the right-hand side of \eqn{tilde N e} is $|J_{e}|$-dimensional Brownian motion because $\lambda_{ei} = 0$ for $i \in J \setminus J_{e}$. Similarly, we have
\begin{align}
\label{eqn:tilde N s}
  \widehat{\vc{N}}_{s}^{(n)}(\cdot) \overset{w}{\longrightarrow} \,\diag(\{\lambda^{\frac 32}_{si} \sigma_{si}\}) \vc{B}_{s}^{*}(\cdot),
\end{align}
where $\vc{B}_{s}^{*}(\cdot)$ is the $d$-dimensional standard Brownian motion. By the modeling assumptions, we can choose $\vc{B}_{s}^{*}(\cdot)$ to be independent of $\vc{B}_{e}^{*}(\cdot)$.

Let $\vc{\chi}^{(n)}(t)$ be the random vector whose $i$-th component is given by
\begin{align*}
  \chi^{(n)}_{i}(t) \equiv \sum_{j \in J} \frac 1{\sqrt{n}} \big(N^{(n)}_{sj}(nt) - \lambda^{(n)}_{sj} nt \big) (p_{ji} - \delta_{ji}),
\end{align*}
then, it follows from \eqn{tilde N s} that
\begin{align}
\label{eqn:tilde chi s}
  \vc{\chi}^{(n)}(\cdot) \overset{w}{\longrightarrow} \Lambda_{\chi} \vc{B}_{s}^{*}(\cdot),
\end{align}
where $\Lambda_{\chi} \Lambda_{\chi}^{\rs{t}}$ is the covariance matrix whose $ik$-entry is given by
\begin{align*}
  \sum_{j \in J} \lambda^{3}_{sj} \sigma_{sj}^{2} (p_{ji} - \delta_{ji}) (p_{jk} - \delta_{jk}).
\end{align*}

We define $\widehat{\vc{\Phi}}^{(n)}(t)$ as the random vector whose $i$-th component is given by
\begin{align*}
  \widehat{\Phi}^{(n)}_{i}(t) \equiv \frac 1{\sqrt{n}} \sum_{j \in J} \big(\Phi_{ji}([nt]) - [nt] p_{ji} \big),
\end{align*}
and define the fluid scaling of $N^{(n)}_{sj}(t)$ as
\begin{align*}
  \ol{N}^{(n)}_{sj}(t) = \frac 1n N^{(n)}_{sj}(nt) .
\end{align*}
Let us compute the covariance matrix of $\widehat{\vc{\Phi}}^{(n)}(1)$, whose $k \ell$-entry is given by
\begin{align}
\label{eqn:CV Phi}
  \frac 1n \dd{E}\Big(\sum_{i, j \in J} & \big(\Phi_{ik}(n) - n p_{ik}\big) \big(\Phi_{j\ell}(n) - n p_{j\ell}\big) \Big) = \sum_{i \in J} p_{ik} (\delta_{k\ell}-p_{i\ell}),
\end{align}
where recall that $\delta_{ij} = 1(i=j)$. Let
\begin{align*}
  \widehat{\vc{\Phi}^{(n)} \circ N^{(n)}_{sj}} (t) = \frac 1{\sqrt{n}} \sum_{j \in J} \big(\Phi_{ji}(N^{(n)}_{sj}(nt)) - N_{sj}^{(n)}(nt) p_{ji} \big).
\end{align*}
Since $\ol{N}^{(n)}_{sj}(t) \to \lambda_{sj} t$ almost surely as $n \to \infty$ and
\begin{align*}
  \frac 1{\sqrt{n}} \sum_{j \in J} \big(\Phi_{ji}(N^{(n)}_{sj}(nt)) - N_{sj}^{(n)}(nt) p_{ji} \big) &= \frac 1{\sqrt{n}} \sum_{j \in J} \big(\Phi_{ji}(n\ol{N}^{(n)}_{sj}(t)) - n \ol{N}_{sj}^{(n)}(t) p_{ji} \big)\\
  & \simeq \frac {1} {\sqrt{n}} \sum_{j \in J} \big(\Phi_{ji}([n \lambda_{sj} t]) - [n \lambda_{sj} t] p_{ji}\big),
\end{align*}
applying \lem{random change} yields
\begin{align}
\label{eqn:tilde Phi}
  \widehat{\vc{\Phi}^{(n)} \circ N^{(n)}_{sj}} (\cdot) \overset{w}{\longrightarrow} \Lambda_{\Phi} \vc{B}_{\Phi}^{*}(\cdot),
\end{align}
where $\vc{B}_{\Phi}^{*}(\cdot)$ is the $d$-dimensional standard Brownian motion which is independent of $\vc{B}_{e}^{*}(\cdot)$ and $\vc{B}_{s}^{*}(\cdot)$, and $\Lambda_{\Phi} \Lambda_{\Phi}^{\rs{t}}$ is the $d$-dimensional covariance matrix whose $ik$-entry is given by
\begin{align*}
  \sum_{j \in J} \lambda_{sj} p_{ji} (\delta_{ik} - p_{jk}),
\end{align*}
which comes from \eqn{CV Phi}.

Finally, the last term in the formula for $\ul{\widehat{X}}^{(n)}_{i}(s)$ is
\begin{align*}
  \sqrt{n} \Big(\lambda^{(n)}_{ei} + \sum_{j \in J} \lambda^{(n)}_{sj} p_{ji} - \lambda^{(n)}_{si} \Big) & = \sqrt{n} \Big(\lambda^{(n)}_{ai} - \sum_{j \in J} \lambda^{(n)}_{aj} p_{ji} + \sum_{j \in J} \lambda^{(n)}_{sj} p_{ji} - \lambda^{(n)}_{si} \Big)\\
 & = \sqrt{n} \sum_{j \in J} (\lambda^{(n)}_{aj} - \lambda^{(n)}_{sj}) (\delta_{ji} - p_{ji})\\
 & = - \left[R \vc{\lambda}_{s} \right]_{i},
\end{align*}
where the last equality comes from \eqn{SR 1}.

Summarizing all the computations, we have that $\widehat{\ul{\vc{X}}}^{(n)}(\cdot)$ weakly converges to the Brownian motion with the covariance matrix $\Sigma$ whose $ik$-entry $\Sigma_{ik}$ is
\begin{align}
\label{eqn:CV 1}
  \Sigma_{ik} = \lambda^{3}_{ei} \sigma_{ei}^{2} \delta_{ik} + \sum_{j \in J} \big(\lambda^{3}_{sj} \sigma_{sj}^{2} (p_{ji} - \delta_{ji}) (p_{jk} - \delta_{jk}) + \lambda_{sj} p_{ji} (\delta_{ik} - p_{jk})\big), \quad i,k \in J,
\end{align}
and drift vector $\vc{\mu}$ whose $i$-entry $\mu_{i}$ is given by
\begin{align}
\label{eqn:drift 1}
  \vc{\mu} = - R \vc{\lambda}_{s}.
\end{align}

Combing this with \eqn{tilde X 1} and \thr{reflection map 1} and applying \thr{continuous}, we arrive at the following theorem.

\begin{theorem}[A version of Theorem 1 of Reiman \cite{Reim1984}] {\rm
\label{thr:GJN 1}
  Assume (\sect{single-class}a)--(\sect{single-class}c), and define the $n$-th GJN by \eqn{Tn 1} and \eqn{Sn 1} and the routing matrix $P$, then the diffusion scaled process $\{\frac 1{\sqrt{n}} \vc{L}^{(n)}(nt); t \ge 0\}$ converges in distribution to the SRBM $\{\vc{Z}(t)\}$ with the primitive data $(\Sigma, - R \vc{\lambda}_{s}, R)$, where $\Sigma$ by \eqn{CV 1}. Thus, this SRBM is always stable.
}\end{theorem}
\begin{remark}
\label{rem:GJN 1}
  If we assume (\sect{diffusion}b) and (\sect{diffusion}c) instead of \eqn{Tn 1} and \eqn{Sn 1}, then the limiting SRBM has primitive data $(\Sigma, - R \vc{c}, R)$, which is stable if and only if $\vc{c} > 0$.
\end{remark}

The following corollary is immediate from Theorems \thrt{BAR} and \thrt{GJN 1} and \eqn{CV 1}.
\begin{corollary}
\label{cor:GJN 1}
  Under the assumptions of \thr{GJN 1}, if $d=1$, then the limiting SRBM has the stationary distribution given by
\begin{align}
\label{eqn:GJN d1}
  \dd{P}( Z_{1}(\infty) \le x) = 1 - \exp\Big( - \frac {2}{\lambda^{2}_{e1}(\sigma^{2}_{e1} + \sigma^{2}_{s1})} x \Big), \qquad x \ge 0,
\end{align}
where $Z_{1}(\infty)$ is a random variable subject to the stationary distribution.
\end{corollary}

Let $\rho^{(n)}_{i} = \lambda^{(n)}_{ei} / \lambda^{(n)}_{si}$. Since \eqn{SR 1} implies that 
\begin{align*}
  \frac 1{\sqrt{n}} = 1 - \rho^{(n)}_{i}, \qquad i \in J,
\end{align*}
we have the following approximation for the stationary queue length process $\vc{L}(\cdot)$.
\begin{align}
\label{eqn:GJN diffusion approx 1}
  \vc{L}(\cdot) \approx \diag(\{1-\rho_{i}\})^{-1} \vc{Z}(\cdot),
\end{align}
where $\vc{Z}(\cdot)$ is the SRBM with the primitive data $(\Sigma, -R\vc{\lambda}_{s}, R)$.

We have not discussed about the workload process $\vc{W}(t)$. Let $\widehat{\vc{W}}^{(n)}(t) = \frac 1{\sqrt{n}} \vc{W}^{(n)}(n t)$, then one may guess that
\begin{align}
\label{eqn:W 2}
  \widehat{\vc{W}}^{(n)}(\cdot) \overset{w}{\longrightarrow} \diag(\{\lambda_{si}\})^{-1} \vc{Z}(\cdot).
\end{align}
and therefore, as $\lambda_{si} -\lambda_{ei} \downarrow 0$ for all $i \in J$,
\begin{align}
\label{eqn:GJN diffusion approx 2}
  \vc{W}(\cdot) \approx \diag(\{\lambda_{si} -\lambda_{ei}\})^{-1} \vc{Z}(\cdot).
\end{align}

We check that \eqn{W 2} indeed holds true. For simplicity, we assume that $\vc{W}^{(n)}(0) = \vc{0}$ and $\vc{L}^{(n)}(0) = \vc{0}$. Then, using \eqn{L 1} for the $n$-th GJN and the fact that
\begin{align*}
  U^{(n)}_{i}(nt) = V^{(n)}_{si}(N^{(n)}_{si}(U^{(n)}_{i}(nt))) + A^{(n)}_{si}(nt),
\end{align*}
where $A^{(n)}_{si}(t)$ is the attained service time of a customer being served at node $i$ if any while $A^{(n)}_{si}(t) = 0$ if no customer in service there, the workload representation \eqn{WL 1} for the $n$-th GJN is rewritten as
\begin{align*}
  \widehat{W}_{i}^{(n)}(t)  & = \frac 1{\sqrt{n}} V^{(n)}_{si}\big(L^{(n)}_{i}(nt) + N^{(n)}_{si}(U^{(n)}_{i}(nt)) \big) - \frac {1}{\sqrt{n}} \big( V^{(n)}_{si}(N^{(n)}_{si}(U^{(n)}_{i}(nt))) + A^{(n)}_{si}(nt) \big) \hspace{2ex}\\
  & = \frac {L^{(n)}_{i}(nt)}{\sqrt{n}} \frac 1{L^{(n)}_{i}(nt)} \sum_{\ell=N^{(n)}_{si}(U^{(n)}_{i}(nt))+1}^{N^{(n)}_{si}(U^{(n)}_{i}(nt)) + L^{(n)}_{i}(nt)} T^{(n)}_{si} - \frac {1}{\sqrt{n}} A^{(n)}_{si}(nt).
\end{align*}
Since $L^{(n)}_{i}(nt) \overset{a.s.}{\longrightarrow} \infty$, which can be checked by fluid approximation, we have
\begin{align*}
  \frac 1{L^{(n)}_{i}(nt)} \sum_{\ell=N^{(n)}_{si}(U^{(n)}_{i}(nt))+1}^{N^{(n)}_{si}(U^{(n)}_{i}(nt)) + L^{(n)}_{i}(nt)} T^{(n)}_{si} \overset{a.s.}{\longrightarrow} \frac 1{\lambda_{si}}, \qquad (n \to \infty).
\end{align*}
Hence, $\frac {1}{\sqrt{n}} A^{(n)}_{si}(nt) \overset{a.s.}{\longrightarrow} 0$ and \thr{GJN 1} yields \eqn{W 2}.

We have derived \thr{GJN 1} along the line due to Reiman \cite{Reim1984}. Some years later, Chen and Mandelbaum \cite{ChenMand1991b} simplified a part of the proof using Skorohod representation of \thr{Skorohod}. This does not have so much benefit for proving \thr{GJN 1}, but it does for extensions or modifications of the GJN. A closed GJN in \cite{ChenMand1991b} is such an example.

\subsection{Extension to many server nodes}
\label{sect:extension to many servers}

\thr{GJN 1} can be extended for the case where each node may have multiple servers and their service times may have different distributions, while keeping $i.i.d$ service times at each server. This extension is already mentioned in \cite{Reim1984}, but was firstly proved in \cite{ChenShan1994}. The basic idea for this extension is called a sandwich method (see also \cite{ChenYe2011} for a single node queue). We present this extension following Chen and Shanthikumar \cite{ChenShan1994}.

Assume that node $i$ has $c_{i}$ servers and the $j$-th server there has service times $\{T_{sij}(\ell)\}$ with finite mean service time $\mu_{sij}^{-1}$ and variance $\sigma_{sij}^{2}$. Define the counting process of departures from this server by $N_{sij}(t)$ as
\begin{align*}
  N_{sij}(t) = \sup\Big\{n \ge 0; \sum_{\ell=1}^{n} T_{sij}(\ell) \le t \Big\},
\end{align*}
and let $U_{ij}(t)$ be the total time when the $j$-th server at node $i$ is busy in the time interval $[0,t]$. Then, for the GJN with multiple servers at nodes, $X_{i}(t)$ of \eqn{X 1} are changed to
\begin{align}
\label{eqn:X 4} 
  X_{i}(t) = & N_{ei}(t) - \lambda_{ei}t + \sum_{j \in J} \sum_{k = 1}^{c_{j}} \big(N_{sjk}(U_{jk}(t)) - \lambda_{sjk} U_{jk}(t) \big)  p_{ji} \nonumber\\
  &  +\sum_{j \in J} \sum_{k = 1}^{c_{j}} \big(\Phi_{ji}(N_{sjk}(U_{jk}(t))) - N_{sjk}(U_{jk}(t)) p_{ji} \big) \nonumber \\
  & - \sum_{k = 1}^{c_{i}} (N_{sik}(U_{ik}(t)) - \lambda_{sik} U_{ik}(t)) + \Big(\lambda_{ei} + \sum_{j \in J} \sum_{k=1}^{c_{j}} \lambda_{sjk} p_{ji} -\sum_{k=1}^{c_{i}}  \lambda_{sik} \Big) t ,
\end{align}
but we can keep $\vc{L}(t)$ of \eqn{L 2} with $\vc{X}(t)$ whose $i$-th component is this $X_{i}(t)$. However, we can not uniquely determine $\vc{L}(t)$ and $\vc{Y}(t)$ by \eqn{L 2} and \eqn{Y 1} obviously because the information for $U_{ik}(t)$ is not sufficient, where
$$Y_{i}(t) = \sum_{k=1}^{c_{i}} \lambda_{sik} (t - U_{ik}(t)).$$

We modify $\vc{X}(t)$ similar to $\ul{\vc{X}}(t)$ of \eqn{X 3} as
\begin{align}
\label{eqn:X 5}
  \ul{X}_{i}(t) = & N_{ei}(t) - \lambda_{ei}t + \sum_{j \in J} \sum_{k = 1}^{c_{j}} \big(N_{sjk}(t) - \lambda_{sjk} t\big) (p_{ji} - \delta_{ij}) \nonumber \\
  & + \sum_{j \in J} \sum_{k = 1}^{c_{j}} \big(\Phi_{ji}(N_{sjk}(t)) - N_{sjk}(t) p_{ji} \big) + \Big(\lambda_{ei} + \sum_{j \in J} \sum_{k = 1}^{c_{j}} \lambda_{sjk} p_{ji} - \sum_{k = 1}^{c_{i}} \lambda_{sik} \Big) t.
\end{align}
We next let $(\phi(\vc{X}), \Psi(\vc{X}))$ be the reflection map for \eqn{L 2} and \eqn{Y 1} with $\vc{X}(\cdot)$ of \eqn{X 5}. Since $\Psi(\vc{X})$ is the minimal solution, we have
\begin{align}
\label{eqn:Psi 2a}
  \Psi(\vc{X})(t) \le \vc{Y}(t).
\end{align}
On the other hand, for the $d$-dimensional vector $\vc{c}$ whose $i$-th component is the number of servers $c_{i}$ at node $i$, it can be shown that $\Psi(\vc{X} - \vc{c})$ is the maximal element in the set:
\begin{align*}
  \big\{ \vc{y} \in \sr{T}; (x_{i}(t) + [R \vc{y}(t)]_{i}) d y_{i}(t) \le c_{i}, i \in J, \forall t \ge 0 \big\},
\end{align*}
where where $\sr{T}$ is the set of non-decreasing functions, and therefore
\begin{align}
\label{eqn:Psi 2b}
  \vc{Y}(t) \le \Psi(\vc{X} - \vc{c})(t).
\end{align}
Thus, we have
\begin{align*}
  \Psi(\vc{X})(t) \le \vc{Y}(t) \le \Psi(\vc{X}-\vc{c})(t).
\end{align*}

We are now ready to introduce the sequence of the GJN with multiple servers at each node. For the $n$-th GJN, we use the following notation at each node $i$. The mean arrival rate and the variance of the inter-arrival time are denoted by $\lambda^{(n)}_{ei}$ and $(\sigma^{(n)}_{ei})^{2}$, which are the same as those for the single server case. For the $j$-th server of node $i$, the mean service rate and the variance of the service time are denoted by $\lambda^{(n)}_{sij}$ and $(\sigma^{(n)}_{sij})^{2}$. We assume that there exist constants $\lambda_{ei}$, $\lambda_{sik}$, $h_{ei}$ and $h_{sik}$ such that
\begin{align}
\label{eqn:root n 1}
 & \lim_{n \to \infty} \sqrt{n} (\lambda^{(n)}_{ei} - \lambda_{ei}) = h_{ei}, \qquad \lim_{n \to \infty} \sqrt{n} (\lambda^{(n)}_{sij} - \lambda_{sij}) = h_{sij},\\
\label{eqn:heavy traffic 1}
 & \lambda_{ei} + \sum_{j \in J} \sum_{k=1}^{c_{j}} \lambda_{sjk} p_{ji} - \sum_{k=1}^{c_{i}} \lambda_{sik} = 0, \qquad \lim_{n \to \infty} \sigma^{(n)}_{ei} = \sigma_{ei}, \qquad \lim_{n \to \infty} \sigma^{(n)}_{sij} = \sigma_{sij},\\
\label{eqn:Lindeberg 4}
 & \sup_{n \ge 1} \dd{E}\Big( (T^{(n)}_{ei})^{2+\delta} + \sum_{j=1}^{c_{j}} (T^{(n)}_{sij})^{2+\delta} \Big) < \infty, \qquad \exists \delta > 0.
\end{align}
Then, \eqn{Psi 2a}, \eqn{Psi 2b} and \eqn{X 5} for the $n$-th GJN concludes the following theorem.

\begin{theorem}[Theorem 5.2 of \cite{ChenShan1994}]
\label{thr:diffusion 1}
Under the assumptions \eqn{root n 1} \eqn{heavy traffic 1} and \eqn{Lindeberg 4}, the diffusion scaled process $\widehat{L}^{(n)}(\cdot) \equiv \{\frac 1{\sqrt{n}} \vc{L}^{(n)}(nt); t \ge 0\}$ converges in distribution to the SRBM $\vc{Z}(\cdot)$ with the primitive data $(\Sigma, \vc{\mu}, R)$, where $\Sigma$ and $\vc{\mu}$ are given by
\begin{align*}
 & \Sigma_{ik} = \lambda^{3}_{ei} \sigma_{ei}^{2} \delta_{ik} + \sum_{j \in J} \sum_{\ell=1}^{c_{j}} \big(\lambda^{3}_{sj\ell} \sigma_{sj\ell}^{2} (p_{ji} - \delta_{ji}) (p_{jk} - \delta_{jk}) + \lambda_{sj\ell} p_{ji} (\delta_{ik} - p_{jk})\big),\\
 & \mu_{i} = h_{ei} + \sum_{j \in J} \sum_{k=1}^{c_{j}} h_{sjk} p_{ji} - \sum_{k=1}^{c_{i}} h_{sik}.
\end{align*}
\end{theorem}

This theorem tells us that the diffusion approximation for the GJN with heterogeneous many servers at all nodes is identical with that for the GNJ with single servers at all nodes with the mean service rate and the variance of service times at node $j$, respectively, given by $\sum_{\ell=1}^{c_{j}} \lambda_{sj\ell}$ and
\begin{align*}
  \Big(\sum_{\ell=1}^{c_{j}} \lambda_{sj\ell}\Big)^{-3} \sum_{\ell=1}^{c_{j}} \lambda^{3}_{sj\ell} \sigma_{sj\ell}^{2}.
\end{align*}
In particular, if servers are identical at each node, that is, $\lambda_{sj\ell} = \lambda_{sj}$ and $\sigma_{sj\ell} = \sigma_{sj}$, then the above variance becomes $\sigma_{sj}^{2}/c_{j}^{2}$, and therefore the result is compatible with Kingman's conjecture for the stationary waiting time in the $GI/GI/k$ queue for $k=c_{j}$, which is proved by K\"{o}llerstr\"{o}m \cite{Koll1974}.

\section{Quality of the diffusion approximation}
\label{sect:quality}

We have discussed about the diffusion approximation through a process limit. This may not imply the weak convergence for the stationary distribution. Let $\widehat{\vc{L}}^{(n)}(\infty)$ and $\vc{Z}(\infty)$ be random vectors subject to the stationary distributions of the diffusion scaled process $\widehat{\vc{L}}^{(n)}(\cdot)$ and the SRBM $\vc{Z}(\cdot)$, respectively. 
Then, the following implication, which is referred to as continuity of the stationary distribution, is important for application because the stationary distribution is widely used for performance evaluation.
\begin{align}
\label{eqn:continuity 1}
  \widehat{\vc{L}}^{(n)}(\cdot) \overset{w}{\longrightarrow} \vc{Z}(\cdot) \quad \mbox{ implies } \quad\widehat{\vc{L}}^{(n)}(\infty) \overset{w}{\longrightarrow} \vc{Z}(\infty).
\end{align}
However, this continuity does not automatically hold in general, and it requires a proof. Similar continuity also may be considered for characteristics of the stationary distribution such as moments and tail asymptotics. Those continuity can be considered as quality support for the diffusion approximation, and therefore their verifications are important for application. In this section, we review the current status for those continuity.

\subsection{Continuity of the stationary distribution and its moments}
\label{sect:continuity distribution}

As already mentioned in \sectn{introduction}, this continuity has been verified by Gamarnik and Zeevi \cite{GamaZeev2006} and Budhiraja and Lee \cite{BudhLee2009}. They use the fact that the tightness for the sequence of $\widehat{\vc{L}}^{(n)}(\infty)$, that is, the sequence of the stationary distributions, is sufficient to prove \eqn{continuity 1}. This is because we can construct a stationary process for the $n$-th GJN taking its stationary distribution as the initial distribution at time $0$ and choose the subsequence of those stationary distributions which weakly converges to a distribution by the tightness.

In \cite{GamaZeev2006}, uniformly light tails over $n$ (see (1) and (2) of \cite{GamaZeev2006} for their precise definitions) are assumed for sequences of the inter-arrival time $T^{(n)}_{ei}$ and the service time $T^{(n)}_{si}$, where we dropped the index $(\ell)$ from $T^{(n)}_{ei}(\ell)$ and $T^{(n)}_{si}(\ell)$. This convention for indexed random variables will be used in what follows if their distributions are only concerned. Then, it is proved that there are positive constants $a_{1}$ and $a_{2}$ for a positive vector $\vc{w} \in \dd{R}^{d}$ such that 
\begin{align*}
  \dd{P}( \br{\vc{w}, \widehat{\vc{L}}^{(n)}(\infty)} > x) \le a_{1} e^{-a_{2}x}, \qquad \forall x > 0.
\end{align*}
See Theorem 7 of \cite{GamaZeev2006}. This is a stronger result more than what is needed for the tightness, and gives an exponential upper bound for the tails of the stationary distribution of the corresponding SRBM, which implies that the moment with any order has the continuity. Furthermore, a finer upper bound is obtained in \cite{GamaZeev2006}. However, those bounds may not be so interesting because a stable SRBM is known to have the light tailed stationary distribution if reflection matrix $R$ is an $\sr{M}$-matrix. 

In \cite{BudhLee2009}, \eqn{continuity 1} is proved under the standard uniformly integrable assumptions for $(T^{(n)}_{ei})^{2}$ and $(T^{(n)}_{si})^{2}$. Furthermore, they proved that, if $\dd{E}([T^{(n)}_{ei}]^{p})$ for $i \in J_{e}$ and $\dd{E}([T^{(n)}_{si}]^{p})$ for $i \in J$ are uniformly finite over $n$ and if certain asymptotic condition on $N_{ei}^{(n)}(\cdot)$ and $N_{si}^{(n)}(\cdot)$ are satisfied (see (i) in (A8.p) of \cite{BudhLee2009}), then the moment with order less than $p$ of the stationary distribution also has the continuity (see Theorem 3.1 of \cite{BudhLee2009}). However, it is known from the finite moment condition for the $GI/GI/1$ queue that $\dd{E}([L_{i}]^{p+1}) = \infty$ if $\dd{E}([T_{si}]^{p}) = \infty$ at least for $i \in J_{e}$ (e.g., see Theorem 2.1 of Section X of \cite{Asmu2003}). Thus, we can not have the continuity for the $(p+1)$-th moment of the stationary distribution in this case because the limiting distribution, which is the stationary distribution of the SRBM, always has light tail and therefore has finite moments for all order (e.g., see \cite{Blan2012}). We may say that the diffusion approximation by SRBM is not good if service time distributions have heavy tails although the continuity holds for the stationary distribution as long as the uniformly integrability holds.

As we demonstrated in \sectn{extension to queue length}, the continuity of the stationary distribution has been directly verified for the single node case of the GJN. As suggested there, the same verification may be possible for a general GJN. Since no time scaling is required for this limiting operation, this continuity has a stronger aspect than the obtained in \cite{BudhLee2009} and \cite{GamaZeev2006}. Furthermore, more information on the convergence may be obtained. In fact, we will see such an example for the tail decay rate of the stationary distribution in the next section. However, a large part of this approach is still open for future study.

\subsection{Continuity of the tail decay rate of the stationary distribution}
\label{sect:continuity tail}

The continuity of the tail decay rate of the stationary distribution is studied for a feedforward queueing network with deterministic service times by Majewski \cite{Maje1998}. This network is assumed to have multi-class customers and FCFS discipline. See \sectn{multi-class} for a multi-class queueing network. It seems that the feedforward assumption and deterministic service times greatly simplifies arguments but its application is limited. It is claimed that the continuity holds for the tail decay rate, precisely, for the large deviations rate function (see Theorem 9 of \cite{Maje1998}).

For a general routing but single-class and two node GJN, the author \cite{Miya2014} recently obtained the tail decay rate of the marginal stationary distribution in an arbitrary direction, assuming a Markov modulated Poisson arrivals and phase type service time distributions, which includes the special case that an arrival process is renewal with a phase type inter-arrival distribution. For simplicity, we assume this renewal arrival.

In the rest of this subsection, we will discuss the continuity on those decay rates. Assume that $d=2$. Recall that $N_{ei}(t)$ is the number of exogenous arrivals at node $i$ during the time period $[0,t]$, and define the logarithmic generating functions similar to \eqn{1-node gamma} as
\begin{eqnarray}
\label{eqn:gamma e}
  \gamma_{ei}(\theta) = \lim_{t \to \infty} \frac 1t \log \dd{E}(e^{\theta N_{ei}(t)}), \qquad i=1,2.
\end{eqnarray}
Because $T_{ei}$ has a phase-type distribution, the assumption (\sect{single-class}d) is satisfied for $T_{ei}$, and therefore these functions are well defined for all $\theta \in \dd{R}$. We similarly define $\gamma_{si}(\theta)$ for the counting process $N_{si}(t)$ for service completions, that is, departures, from node $i$ during a busy period. 

Similarly in \sectn{extension to queue length}, we denote the moment generating functions of $T_{ui}$ is denoted by
\begin{align*}
 \widetilde{F}_{ui}(\theta) = \dd{E}(e^{\theta T_{ui}}), \qquad u=e,s, \; i=1,2,
\end{align*}
then we have the following facts similar to \lem{1-node logarithmic}.
\begin{lemma}
\label{lem:logarithmic 1}
  $\widetilde{F}_{ui}$ has a unique inverse $\widetilde{F}_{ui}^{-1}$ whose domain is $\dd{R}$, and, for $u=e,s$ and $i=1,2$,
\begin{align}
\label{eqn:logarithmic 1}
  \gamma_{ui}(\theta) = - \widetilde{F}_{ui}^{-1}(e^{-\theta}), \qquad \theta \in \dd{R}.
\end{align}
\end{lemma}

We next consider a logarithmic moment generating function for customers' movement at their departure. Recall that departures from node $i$ are routed to node $j$ according to the Markovian routing function $\Phi_{ij}(\cdot)$. Hence, if a server at node $i$ is busy for all time up to $t$, then we define the logarithmic moment generating function for departing and routing customers for busy node $i$ as
\begin{eqnarray*}
  \gamma_{di}(\vc{\theta}) = \lim_{t \to \infty} \frac 1t \log \dd{E}(e^{-\theta_{i} N_{si}(t) + \theta_{3-i} \Phi_{i(3-i)}(N_{si}(t))}), \qquad \vc{\theta} \equiv (\theta_{1}, \theta_{2}) \in \dd{R}^{2}.
\end{eqnarray*}
Let $t_{i}(\vc{\theta}) = e^{-\theta_{i}} (p_{i0} + p_{ij} e^{\theta_{j}})$ for $j = 3-i$, then it is not hard to see that
\begin{align}
\label{eqn:gamma d}
  \gamma_{di}(\vc{\theta}) &= \lim_{t \to \infty} \frac 1t \log \dd{E}(t_{i}(\vc{\theta})^{N_{si}(t)}) \nonumber \\
  &= \lim_{t \to \infty} \frac 1t \log \dd{E}(e^{N_{si}(t) \log t_{i}(\vc{\theta})}) = \gamma_{si}(\log t_{i}(\vc{\theta})).
\end{align}

We are now ready to present main results. For this, we need some notations. Let
\begin{align*}
 & \gamma_{+}(\vc{\theta}) = \gamma_{e1}(\theta_{1}) + \gamma_{e2}(\theta_{2}) + \gamma_{d1}(\vc{\theta}) + \gamma_{d2}(\vc{\theta}),
\end{align*}
which stands for the logarithmic generating function for the total net flow. The following notations are geometric objects.
\begin{align*}
 & \Gamma = \{\vc{\theta} \in \dd{R}^{2}; \gamma_{+}(\vc{\theta}) < 0 \}, \qquad \partial \Gamma = \{\vc{\theta} \in \dd{R}^{2}; \gamma_{+}(\vc{\theta}) = 0 \},,\\
 & \Gamma_{i} = \{\vc{\theta} \in \Gamma; t_{3-i}(\vc{\theta}) > 1 \}, \qquad \partial \Gamma_{i} = \{\vc{\theta} \in \partial \Gamma; t_{3-i}(\vc{\theta}) > 1 \}, \qquad i=1,2,\\
 & \Gamma_{\max} = \{\vc{\theta} \in \dd{R}^{2}; \exists \vc{\theta}' \in \partial \Gamma, \vc{\theta} < \vc{\theta}'\}.
\end{align*}

Let $\varphi(\vc{\theta}) \equiv \dd{E}(e^{\br{\vc{\theta}, \vc{L}(\infty)}})$ be the moment generating function of $\vc{L}(\infty)$ subject to the marginal stationary distribution with respect to the numbers of customers at nodes 1 and 2, and define $\sr{D}$ as
\begin{eqnarray*}
  \sr{D} = \mbox{the interior of } \{\vc{\theta} \in \dd{R}^{2}; \varphi(\vc{\theta}) < \infty\},
\end{eqnarray*}
which is referred to as a convergence domain of $\varphi$. Denote the two extreme points of $\Gamma_{i}$ by
\begin{eqnarray*}
  && \vc{\theta}^{(i,\cp)} = \arg_{\vc{\theta} \in \dd{R}^{2}} \sup\{ \theta_{i} \ge 0; \vc{\theta} \in \partial \Gamma_{i} \} , \quad i=1,2.
\end{eqnarray*}
Using these points, we define the vector $\vc{\tau}$ by
\begin{eqnarray*}
 && \tau_{1} = \sup\{\theta_{1} \in \dd{R}; \vc{\theta} \in \partial \Gamma_{1}; \theta_{2} < \theta^{(2,\cp)}_{2} \},\\
 && \tau_{2} = \sup\{\theta_{2} \in \dd{R}; \vc{\theta} \in \partial \Gamma_{2}; \theta_{1} < \theta^{(1,\cp)}_{1} \}.
\end{eqnarray*}

We refer to the following results (see \cite{Miya2014} for their proofs).

\begin{theorem}[Theorems 3.1 and 3.3 of Miyazawa \cite{Miya2014}]
\label{thr:decay rate 1}
  For the two node generalized Jackson network under phase type setting, if it is stable, then $ \sr{D} = \{ \vc{\theta} \in \Gamma_{\max}; \vc{\theta} < \vc{\tau} \}$, and, for any non-zero vector $\vc{c} \ge \vc{0}$,
\begin{eqnarray}
\label{eqn:decay rate 1}
 && \lim_{x \to \infty} \frac 1x \log \dd{P}(\br{\vc{c}, \vc{L}(\infty)} > x) = - \sup \{u > 0; u \vc{c} \in \sr{D} \}.
\end{eqnarray}
\end{theorem}

Note that this theorem shows that the decay rates are determined by \eqn{gamma e} and \eqn{gamma d}. Hence, we can expect that \thr{decay rate 1} hold beyond the phase-type setting and more general routing mechanism.

It is also notable that these decay rates are obtained in the exactly same way as those of the stationary distribution of the two-dimensional SRBM (see Theorems 2.1, 2.2 and 2.3 of \cite{DaiMiya2011}). To prove the continuity of the decay rates, we here briefly present results in \cite{DaiMiya2011}.

Let $\{(\vc{Z}(t), \vc{Y}(t))\}$ be the two dimensional SRBM with the primitive data $(\Sigma, \vc{\mu}, R)$, and assume that $R$ is a completely-$\sr{S}$ matrix and the stability condition $R^{-1} \vc{\mu} < \vc{0}$ is satisfied. Define the following functions and sets.
\begin{align*}
 & \widetilde{\gamma}_{+}(\vc{\theta}) = \frac 12 \br{\vc{\theta}, \Sigma \vc{\theta}} + \br{\vc{\mu}, \vc{\theta}}, \qquad \widetilde{\gamma}_{i}(\vc{\theta}) = \br{R^{[i]}, \vc{\theta}}, \qquad i = 1,2,\\
 & \widetilde{\Gamma} = \{\vc{\theta} \in \dd{R}^{2}; \widetilde{\gamma}_{+}(\vc{\theta}) < 0\}, \qquad \widetilde{\Gamma}_{i} = \{\vc{\theta} \in \widetilde{\Gamma}; \widetilde{\gamma}_{3-i}(\vc{\theta}) < 0\},\\
 & \widetilde{\Gamma}_{\max} = \{\vc{\theta} \in \dd{R}^{2}; \exists \vc{\theta}' \in \widetilde{\Gamma}, \vc{\theta} < \vc{\theta}'\}, 
\end{align*}
where we recall that $R^{[i]}$ is the $i$-th column of $R$. Let $\vc{Z}(\infty)$ be a random vector subject to the stationary distribution of $\{\vc{Z}(t)\}$. Note that this stationary distribution is uniquely determined by $\widetilde{\gamma}_{+}$ and $\widetilde{\gamma}_{i}$ (see \rem{BAR} for its reason). Recall that $\tilde{\varphi}(\vc{\theta}) = \dd{E}(e^{\br{\vc{\theta}, \vc{Z}(\infty)}})$. Define the convergence domain of $\tilde{\varphi}(\vc{\theta})$ as
\begin{align*}
  \widetilde{\sr{D}} = \mbox{the interior of } \{\vc{\theta} \in \dd{R}^{2}; \tilde{\varphi}(\vc{\theta}) < \infty\}.
\end{align*}
We also define $\widetilde{\vc{\theta}}^{(i,\Gamma)}$ and $\widetilde{\vc{\tau}}$ in the exactly same way as $\vc{\theta}^{(i,\Gamma)}$ and $\vc{\tau}$. Then, similar to \thr{decay rate 1}, we have the following results from Theorems 2.1, 2.2 and 2.3 in \cite{DaiMiya2011}, in which finer tail asymptotics are obtained, but we do not use them here.

\begin{theorem}[Dai and Miyazawa \cite{DaiMiya2011}]
\label{thr:decay rate 2}
  For the two node generalized Jackson network under phase type setting, if it is stable, then $ \widetilde{\sr{D}} = \{ \vc{\theta} \in \widetilde{\Gamma}_{\max}; \vc{\theta} < \widetilde{\vc{\tau}} \}$, and, for any non-zero vector $\vc{c} \ge \vc{0}$,
\begin{eqnarray}
\label{eqn:decay rate 2}
 && \lim_{x \to \infty} \frac 1x \log \dd{P}(\br{\vc{c}, \vc{Z}(\infty)} > x) = - \sup \{u > 0; u \vc{c} \in \widetilde{\sr{D}} \}.
\end{eqnarray}
\end{theorem}

We now consider the sequence of the GJNs whose joint queue length processes under diffusion scaling weakly converge to a stable SRBM. For this, we use the same notations and assumptions in \sectn{single-class}. For example, recall that $\lambda^{(n)}_{ei}, \lambda^{(n)}_{si}$ and $(\sigma^{(n)}_{ei})^{2}, (\sigma^{(n)}_{si})^{2}$ are the arrival and service rates and the variances of the inter-arrival and service times, respectively. We assume the heavy traffic condition (\sect{diffusion}b) with $c_{i} = \lambda_{ai} - \lambda_{si}$ and regularity condition (\sect{diffusion}c). Then, the covariance matrix $\Sigma$ and the drift vector $\vc{\mu}$ of this SRBM are given by \eqn{CV 1} and \eqn{drift 1}, respectively, and $R = I - P^{\rs{t}}$. The logarithmic generating function of the $n$-th system under diffusion scaling is defined as
\begin{align}
\label{eqn:logarithmic diffusion}
  \widehat{\gamma}^{(n)}_{ui}(\theta) = \lim_{t \to \infty} \frac 1{t} \log \dd{E}(e^{\theta \frac 1{\sqrt{n}} N^{(n)}_{ui}(nt)}) = n \gamma^{(n)}_{ui}\Big(\frac 1{\sqrt{n}} \theta\Big), \qquad u = e, s.
\end{align}
Hence, by Lemmas \lemt{Taylor 1} and \lemt{logarithmic 1}, for $u = e,s$,
\begin{align*}
  \widehat{\gamma}^{(n)}_{ui}(\theta) = - n (\widetilde{F}^{(n)}_{ui})^{-1}\big(e^{-\frac {\theta}{\sqrt{n}}}\big) =  \sqrt{n}\lambda_{ui}^{(n)} \theta + \frac 12 (\lambda^{(n)}_{ui})^{3} (\sigma^{(n)}_{ui})^{2} \theta^{2} + o(1),
\end{align*}
and therefore it follows from \eqn{gamma d} that, for $i=1,2$ and $j = 3 - i$, 
\begin{align*}
  \widehat{\gamma}^{(n)}_{di}(\vc{\theta}) &= - n (\widetilde{F}^{(n)}_{si})^{-1}\Big(e^{-\log t_{i}\left(\frac {1}{\sqrt{n}}\vc{\theta}\right)}\Big)\\
  & = n \left(\lambda_{si}^{(n)} \log t_{i}\left(\frac {1}{\sqrt{n}}\vc{\theta}\right) + \frac 12 (\lambda^{(n)}_{si})^{3} (\sigma^{(n)}_{si})^{2} \log^{2} t_{i}\left(\frac {1}{\sqrt{n}}\vc{\theta}\right) \right) + o(1)\\
  & = \sqrt{n} \lambda^{(n)}_{si} (-\theta_{i} + p_{ij} \theta_{j}) + \frac 12 \lambda^{(n)}_{si} p_{ij} (1-p_{ij}) \theta_{j}^{2} + \frac 12 (\lambda^{(n)}_{si})^{3}(\sigma^{(n)}_{si})^{2} (\theta_{i} - \theta_{j})^{2} + o(1).
\end{align*}

From these facts, we can see that
\begin{align}
\label{eqn:gamma continuity}
  \lim_{n \to \infty} \widehat{\gamma}^{(n)}_{+}(\vc{\theta}) = \widetilde{\gamma}_{+}(\vc{\theta}), \qquad \lim_{n \to \infty} (- (\lambda^{(n)}_{si} \sqrt{n})^{-1}) \widehat{\gamma}^{(n)}_{di}(\vc{\theta}) = \widetilde{\gamma}_{i}(\vc{\theta}), \quad i = 1,2.
\end{align}
Thus, in the view of Theorems \thrt{decay rate 1} and \thrt{decay rate 2}, we are very close to have the continuity of the decay rates. However, the second equation of \eqn{gamma continuity} is a bit short to prove the following theorem. We give a complete proof in \app{decay}.

\begin{theorem}
\label{thr:decay rate continuity}
  For a sequence of two node GJNs with phase type $T^{(n)}_{e} \equiv T_{e}$ and $T^{(n)}_{s}$, if they are stable and satisfy the condition \eqn{Sn 1}, then the convergence domain $\sr{D}^{(n)}$ of the moment generating function of the stationary distribution of the $n$-th GJN converges to the domain $\widetilde{\sr{D}}$ for the limiting SRBM, and therefore the decay rate of $\br{\vc{c}, \vc{L}^{(n)}(\infty)}$ for any non-zero $\vc{c} \ge 0$ converges to that of the limiting SRBM.
\end{theorem}

%------------------------------------------------------ 
\section{Multi-class queueing network}
\label{sect:multi-class}
%------------------------------------------------------ 

We discuss how diffusion approximation is obtained for a multi-class queueing network with FCFS service discipline, following Williams \cite{Will1998a}. We modify the GJN for this. We denote the set of customer types by $K \equiv \{1,2,\ldots,c\}$. Type $k$ customers are served only at node $\beta(k)$, and changed to type $k'$ customer after service completion with probability $p_{kk'}$ independent of everything else. Similarly, $T_{ei}$ and $T_{si}$ are changed to $T_{ek}$ and $T_{sk}$, respectively, for $k \in K$.

Let $N_{ek}(t), N_{ak}(t), N_{sk}(t), N_{dk}(t)$ be the counting processes for exogenous arrivals, total arrivals, potential service completions, departures of type $k$ customers. Let $L_{k}(t), U_{k}(t)$ be the number of customers and the accumulated service time for type $k$ customers. Let $W_{i}(t), Y_{i}(t)$ be the workload and the accumulated idle time at node $i$ at time $t$, respectively. For simplicity, we assume that $\vc{L}(0) = \vc{0}$ and $\vc{W}(0) = \vc{0}$. Then, they satisfy the following set of equations for $i \in J$ and $k \in K$. Recall that the composition $f \circ g$ is defined as $f \circ g (x) = f(g(x))$ for appropriate functions $f$ and $g$.
\begin{align}
\label{eqn:mul Na 1}
 & N_{ak}(t) = N_{ek}(t) + \sum_{k' \in K} (\Phi_{k'k} \circ N_{dk'})(t),\\
\label{eqn:mul L 1}
 & L_{k}(t) = N_{ak}(t) - N_{dk}(t),\\
\label{eqn:mul Nd 1}
 & N_{dk}(t) = (N_{sk} \circ U_{k})(t),\\
\label{eqn:mul W 1}
 & W_{i}(t) = \sum_{k \in K} c_{ik} (V_{sk} \circ N_{ak})(t) - t + Y_{i}(t),\\
\label{eqn:mul Y 1}
 & Y_{i}(t) = t - \sum_{k \in K} c_{ik} U_{k}(t),\\
\label{eqn:mul WY 1}
 & \int_{0}^{t} W_{i}(u) dY_{i}(u) = 0, 
\end{align}
where $c_{ik} = 1(\beta(k)=i)$. One can see that \eqn{mul L 1} with \eqn{mul Na 1}, \eqn{mul W 1} with \eqn{mul Y 1}, and \eqn{mul WY 1} are parallel to \eqn{L 1}, \eqn{WL 1} and \eqn{Y 2} of a single-class queueing network. An essential feature here is that the workload $\vc{W}(t)$ plays a key role while the queue length $\vc{L}(t)$ did so in the single-class case. It also indicates the lack of information because $\vc{L}(t)$ is higher dimensional than $\vc{W}(t)$. To fill this gap, we need to specify service discipline among customer types at each node.

Once the service discipline is specified, we can fully characterize $\vc{W}(t)$ and $\vc{L}(t)$ by \eqn{mul Na 1}--\eqn{mul WY 1}. For example, if the service discipline is FCFS independent of types, then the following conditions are sufficient.
\begin{align}
\label{eqn:mul FCFS}
  N_{ak}(t) = N_{dk}(t + W_{\beta(k)}(t)),
\end{align}
which represents that the type $k$ customers arriving up to time $t$ will depart up to time $t + W_{\beta(k)}(t)$. However, this is said in principle, and it is almost impossible to compute them. On the other hand, one may guess that the number of each type customers in node $i$ would be proportional to their arrival rate in heavy traffic because the queue there would be very long and service does not discriminate types. This is called a state-space collapse, which will be formally defined below. Thus, if the state-space collapse occurs, then we may compute $\vc{L}(t)$ from the lower dimensional $\vc{W}(t)$. This is exactly a key idea to get diffusion approximation for a multi-class queueing network. 

We now introduce the sequence of the multi-class GJN's. They are indexed by natural numbers similar to those of the single-class GJN's. Thus, the $n$-th system have the following notations. $T^{(n)}_{ek}(\ell)$ and $T^{(n)}_{sk}(\ell)$ are the inter-arrival (from the outside) and service times of type $k$ customers. It is assume that they are $i.i.d.$ for each fixed $n$ and $k \in K$. Counting processes of type $k$ customers and their mean rates are denoted by $N^{(n)}_{uk}(t)$ and $\lambda^{(n)}_{uk}$ for $u =e, a, s, d$. Routing function, the number of type $k$ customers, the total workload brought by $\ell$ customers of type $k$ are denoted by $\Phi^{(n)}_{kk'}(\ell)$, $L^{(n)}_{k}(t)$, $V^{(n)}_{sk}(\ell)$ for $k \in K$. The workload and total idle time at node $i$ are denoted by $W^{(n)}_{i}(t)$, $Y^{(n)}_{i}(t)$ for $i \in J$. Vectors and matrix whose entries are given by them are denoted by $\vc{N}^{(n)}_{u}(t)$ and $\vc{\lambda}^{(n)}_{u}$ for $u =e, a, s, d$, $\vc{\Phi}^{(n)}(\ell)$, $\vc{L}^{(n)}(t)$, $\vc{V}^{(n)}_{s}(\ell)$, $\vc{W}^{(n)}(t)$ and $\vc{Y}^{(n)}(t)$, respectively. We also introduce $|K| \times |K|$ matrices $P \equiv \{p_{kk'}; k,k' \in K\}$ and $M^{(n)} \equiv \mbox{diag} (\{m^{(n)}_{sk}; k \in K\})$, $d \times K$ matrix $C \equiv \{1(\beta(k)=i); i \in J, k \in K\}$, where $m^{(n)}_{sk} = (\lambda^{(n)}_{sk})^{-1}$. We define the traffic intensity $\rho^{(n)}_{k} = \lambda^{(n)}_{ak} m^{(n)}_{sk}$, where the total arrival rate of type $k$ customer $\lambda_{ak}$ is given by
\begin{align*}
  \vc{\lambda}^{(n)}_{a} = (I - P^{\rs{t}})^{-1} \vc{\lambda}^{(n)}_{e}.
\end{align*}

We define fluid scaling versions of the above processes as
\begin{align*}
 & \ol{N}^{(n)}_{uk}(t) = \frac 1n N^{(n)}_{uk}(nt), \qquad \ol{\Phi}^{(n)}_{kk'}(t) = \frac 1n \Phi^{(n)}_{kk'}([nt]), \\
 & \ol{L}^{(n)}_{k}(t) = \frac 1n L^{(n)}_{k}(nt), \qquad \ol{V}^{(n)}_{sk}(t) = \frac 1n V^{(n)}_{sk}([nt]),\\
 & \ol{W}^{(n)}_{i}(t) = \frac 1n W^{(n)}_{i}(nt), \qquad \ol{Y}^{(n)}_{i}(t) = \frac 1n Y^{(n)}_{i}(nt).
\end{align*}
Similarly, diffusion scaling versions are defined as
\begin{align*}
 & \widehat{N}^{(n)}_{uk}(t) = \frac 1{\sqrt{n}} \big(N^{(n)}_{uk}(nt) - \lambda^{(n)}_{uk} nt\big), \qquad \widehat{\Phi}^{(n)}_{kk'}(t)= \frac 1{\sqrt{n}} \big(\Phi^{(n)}_{kk'}([nt]) - p_{kk'}[nt]\big) \\
 & \widehat{L}^{(n)}_{k}(t) = \frac 1{\sqrt{n}} L^{(n)}_{k}(nt),\qquad \ol{V}^{(n)}_{sk}(t) = \frac 1{\sqrt{n}} \big( V^{(n)}_{sk}([nt]) - m^{(n)}_{sk} [nt]) \big),\\
 & \widehat{W}^{(n)}_{i}(t) = \frac 1{\sqrt{n}} W^{(n)}_{i}(nt), \qquad \widehat{Y}^{(n)}_{i}(t) = \frac 1{\sqrt{n}} Y^{(n)}_{i}(nt),\\
 & \widehat{V^{(n)}_{sk} \circ N^{(n)}_{ak}}(t) = \frac 1{\sqrt{n}} \big( V^{(n)}_{sk} \circ N^{(n)}_{ak}(nt) - \lambda^{(n)}_{ak} m^{(n)}_{sk} nt\big),\\
 & \widehat{\Phi^{(n)}_{kk'} \circ N^{(n)}_{dk}}(t) = \frac 1{\sqrt{n}} \big( \Phi^{(n)}_{kk'} \circ N^{(n)}_{dk}(nt) - \lambda^{(n)}_{dk} p^{(n)}_{kk'} nt\big).
\end{align*}

From those definitions, we have
\begin{align}
\label{eqn:dif Phi 1}
  \widehat{\Phi}^{(n)}_{k'k}(\ol{N}^{(n)}_{dk'}(t)) &= \frac 1{\sqrt{n}} \big( \Phi^{(n)}_{k'k}(N^{(n)}_{dk'}(nt)) - p_{k'k} N^{(n)}_{dk'}(nt)\big) \nonumber\\
 &= \widehat{\Phi^{(n)}_{k'k} \circ N^{(n)}_{dk'}}(t) - p_{k'k} \widehat{N}^{(n)}_{dk'}(t).
\end{align}
Similarly,
\begin{align}
\label{eqn:dif Vs 1}
  \widehat{V}^{(n)}_{sk}(\ol{N}^{(n)}_{ak}(t)) &= \widehat{V^{(n)}_{sk} \circ N^{(n)}_{ak}}(t) - m^{(n)}_{sk} \widehat{N}^{(n)}_{ak}(t).
\end{align}

Substituting \eqn{dif Phi 1} into the $n$-th diffusion scaling version of \eqn{mul Na 1}, we have
\begin{align*}
  \widehat{N}^{(n)}_{ak}(t) &= \widehat{N}^{(n)}_{ek}(t) + \sum_{k' \in K} \widehat{\Phi^{(n)}_{k'k} \circ N^{(n)}_{dk'}}(t) \nonumber\\
  &= \widehat{N}^{(n)}_{ek}(t) + \sum_{k' \in K} \big(\widehat{\Phi}^{(n)}_{k'k}(\ol{N}^{(n)}_{dk'}(t)) + p_{k'k} \widehat{N}^{(n)}_{dk'}(t)\big).
\end{align*}
Substituting the $n$-th diffusion scaling version of \eqn{mul L 1} into this equation, we have
\begin{align*}
  \widehat{N}^{(n)}_{ak}(t) = \widehat{N}^{(n)}_{ek}(t) + \sum_{k' \in K} \big(\widehat{\Phi}^{(n)}_{k'k}(\ol{N}^{(n)}_{dk'}(t)) + p_{k'k} (\widehat{N}^{(n)}_{ak'}(t) - \widehat{L}^{(n)}_{k'}(t)) \big),
\end{align*}
which is expressed as, in vector-matrix notation,
\begin{align*}
  \widehat{\vc{N}}^{(n)}_{a}(t) = \widehat{\vc{N}}^{(n)}_{e}(t) + \sum_{k' \in K} \widehat{\vc{\Phi}}^{(n)}_{k'}(\ol{N}^{(n)}_{dk'}(t)) + P^{\rs{t}} (\widehat{\vc{N}}^{(n)}_{a}(t) - \widehat{\vc{L}}^{(n)}(t)).
\end{align*}
Hence, multiplying $Q \equiv (I - P^{\rs{t}})^{-1}$ from the left, we have
\begin{align}
\label{eqn:dif Na 1}
  \widehat{\vc{N}}^{(n)}_{a}(t) = Q \widehat{\vc{N}}^{(n)}_{e}(t) + Q \sum_{k' \in K} \widehat{\vc{\Phi}}^{(n)}_{k'}(\ol{N}^{(n)}_{dk'}(t)) - Q P^{\rs{t}} \widehat{\vc{L}}^{(n)}(t).
\end{align}

We further substitute \eqn{dif Na 1} into the $n$-th diffusion version of \eqn{mul W 1}, then, using \eqn{dif Vs 1} and again \eqn{dif Na 1}, arrive at
\begin{align}
\label{eqn:dif W 1}
  \widehat{\vc{W}}^{(n)}(t) &= C \widehat{\vc{V}^{(n)}_{s} \circ \vc{N}^{(n)}_{a}}(t) + \sqrt{n} t (\vc{\rho}^{(n)} - \vc{1}) + \widehat{\vc{Y}}^{(n)}(t) \nonumber\\
  &= C \big(\widehat{\vc{V}}^{(n)}_{s} (\ol{\vc{N}}^{(n)}_{a}(t)) + M^{(n)} \widehat{\vc{N}}^{(n)}_{a}(t) \big) + \sqrt{n} t (\vc{\rho}^{(n)} - \vc{1}) + \widehat{\vc{Y}}^{(n)}(t) \nonumber\\
  & = C \big(\widehat{\vc{V}}^{(n)}_{s} (\ol{\vc{N}}^{(n)}_{a}(t)) + M^{(n)} Q (\widehat{\vc{N}}^{(n)}_{e}(t) + \sum_{k' \in K} \widehat{\vc{\Phi}}^{(n)}_{k'}(\ol{N}^{(n)}_{dk'}(t)) \big) \nonumber\\
  & \hspace{20ex} - CM^{(n)}Q P^{\rs{t}} \widehat{\vc{L}}^{(n)}(t) + \sqrt{n} t (\vc{\rho}^{(n)} - \vc{1}) + \widehat{\vc{Y}}^{(n)}(t).
\end{align}

We now formally define the state space collapse according to \cite{Will1998}.
\begin{definition}
\label{dfn:collapse 1}
  State space collapse is said to hold for $K \times J$ nonnegative matrix $\Delta^{(n)}$, which is called a lifting matrix, if $ki$-entry of $\Delta^{(n)}$ vanishes if $\beta(k) \not= i$ and, for each $t \ge 0$,
\begin{align}
\label{eqn:state space collapse 1}
  \| \widehat{\vc{L}}^{(n)}(\cdot) - \Delta^{(n)} \widehat{\vc{W}}^{(n)}(\cdot) \|_{t} \overset{P}{\longrightarrow} 0, \qquad n \to \infty,
\end{align}
where $\overset{P}{\longrightarrow}$ stands for convergence in probability.
\end{definition}

Assume that state space collapse holds for the sequence of lifting matrices $\Delta^{(n)}$, and $\Delta^{(n)}$ converges to $\Delta \ne 0$. In particular, we assume FCFS service discipline. Then, as we expected, it can be proved that the condition \eqn{state space collapse 1} holds if
\begin{align}
\label{eqn:Delta FCFS}
  \Delta^{(n)}_{ki} = \frac {\lambda^{(n)}_{ak}} {\sum_{\beta(k) = i} \lambda^{(n)}_{ak} m^{(n)}_{sk}},
\end{align}
and therefore
\begin{align}
\label{eqn:CMD 1}
  C M^{(n)} \Delta^{(n)} = I.
\end{align}

We write the condition \eqn{state space collapse 1} as
\begin{align*}
  \widehat{\vc{L}}^{(n)}(t) = \Delta^{(n)} \widehat{\vc{W}}^{(n)}(t) + \widehat{\vc{\epsilon}}^{(n)}(t) \quad \mbox{with} \quad \|\widehat{\vc{\epsilon}}^{(n)}(t)\|_{t} \overset{P}{\longrightarrow} 0, \quad n \to \infty,
\end{align*}
and substitute $\widehat{\vc{L}}^{(n)}(t)$ into \eqn{dif W 1}, then we have
\begin{align}
\label{eqn:dif W 2}
  (I + CM^{(n)}QP^{\rs{t}} \Delta^{(n)}) \widehat{\vc{W}}^{(n)}(t) = & C \big(\widehat{\vc{V}}^{(n)}_{s} (\ol{\vc{N}}^{(n)}_{a}(t)) + M^{(n)} Q (\widehat{\vc{N}}^{(n)}_{e}(t) + \sum_{k' \in K} \widehat{\vc{\Phi}}^{(n)}_{k'}(\ol{N}^{(n)}_{dk'}(t)) \big) \nonumber\\
  & - CM^{(n)} Q P^{\rs{t}} \widehat{\vc{\epsilon}}^{(n)}(t) + \sqrt{n} t (\vc{\rho}^{(n)} - \vc{1}) + \widehat{\vc{Y}}^{(n)}(t).
\end{align}
Define $G^{(n)}$ as
\begin{align*}
  G^{(n)} = CM^{(n)}QP^{\rs{t}} \Delta^{(n)}.
\end{align*}
Because of \eqn{CMD 1}, we have
\begin{align}
\label{eqn:I+Gn}
  I + G^{(n)} = CM^{(n)}(I + QP^{\rs{t}})\Delta^{(n)} = CM^{(n)} Q \Delta^{(n)}.
\end{align}

We now assume the heavy traffic conditions: For some constant $c_{k} \in \dd{R}$,
\begin{align}
\label{eqn:heavy 1}
  \lambda^{(n)}_{ek} \to \lambda_{ek}, \qquad m^{(n)}_{sk} \to m_{sk}, \qquad \sqrt{n} (\rho^{(n)}_{k} - 1) \to c_{k},
\end{align}
where $\rho_{k} = \lambda_{ak} m_{sk}$ and $\lambda_{ak}$ is given by 
\begin{align*}
  \vc{\lambda}_{a} = Q \vc{\lambda}_{e}.
\end{align*}
For diffusion approximation, we assume the finiteness of variances. Let
\begin{align*}
  (\sigma^{(n)}_{ek})^{2} = \dd{E}((T^{(n)}_{ek}(1))^{2}), \qquad (\sigma^{(n)}_{sk})^{2} = \dd{E}((T^{(n)}_{sk}(1))^{2}), \qquad k \in K,
\end{align*}
where we use the convention that $(\sigma^{(n)}_{ek})^{2} = 0$ if $\lambda^{(n)}_{ek} = 0$, and assume that
\begin{align}
\label{eqn:mul variance 1}
  \sigma^{(n)}_{ek} \to \sigma_{ek}, \qquad \sigma^{(n)}_{sk} \to \sigma_{sk}.
\end{align}
Because we did not choose specific sequences for the inter-arrival times $T^{(n)}_{ek}(\ell)$ and the service times $T^{(n)}_{sk}(\ell)$ as we did for the single-class GJN, we here require the following Lindeberg type conditions.
\begin{align}
\label{eqn:Lindeberg 1}
  \max_{k \in K} \sup_{n \in \dd{N}} \dd{E}\big( (T^{(n)}_{ek}(1))^{2}; T^{(n)}_{ek}(1) > n) \to 0, \\
\label{eqn:Lindeberg 2}
  \max_{k \in K} \sup_{n \in \dd{N}} \dd{E}\big( (T^{(n)}_{sk}(1))^{2}; T^{(n)}_{sk}(1) > n) \to 0.
\end{align}

Under these assumptions, it can be proved the following convergence in distribution.
\begin{align}
\label{eqn:dif Vs 2}
 & \widehat{\vc{V}}^{(n)}_{s} (\ol{\vc{N}}^{(n)}_{a}(t)) \overset{w}{\longrightarrow} \diag(\{\lambda^{\frac 12}_{sk} \sigma_{sk}\}) \vc{B}_{s}^{*}(\cdot)\\
\label{eqn:dif Ne 1}
 & \widehat{\vc{N}}^{(n)}_{e}(t) \overset{w}{\longrightarrow} \diag(\{\lambda_{ek}^{\frac 32} \sigma_{ek}\}) \vc{B}_{e}^{*}(\cdot)\\
\label{eqn:dif Phi k 1}
 & \widehat{\vc{\Phi}}^{(n)}_{k'}(\ol{N}^{(n)}_{dk'}(t)) \overset{w}{\longrightarrow} \lambda_{ak'}^{\frac 12} \Gamma_{k'} \vc{B}_{k'}^{*}(\cdot),
\end{align}
where $\vc{B}_{s}^{*}(\cdot)$, $\vc{B}_{e}^{*}(\cdot)$ and $\vc{B}_{k}^{*}(\cdot)$ are independent $|K|$, $|K_{e}|$ and $|K|$ dimensional standard Brownian motions, and $\Gamma_{k'} \Gamma_{k'}^{\rs{t}}$ is the $K \times K$ covariance matrix defined by
\begin{align*}
  [\Gamma_{k'} \Gamma_{k'}^{\rs{t}}]_{\ell \ell'} = \left\{\begin{array}{ll}
  p_{k \ell} (1 - p_{k \ell}), \quad & \ell' = \ell,\\
  - p_{k \ell} p_{k \ell'}, & \ell' \ne \ell.
  \end{array} \right.
\end{align*}

Let $G = \lim_{n \to \infty} G^{(n)}$, and assume that $I+G$ is invertible. Let $R \equiv (I+G)^{-1}$. Taking \eqn{dif W 2} and \eqn{dif Vs 2}--\eqn{dif Phi k 1} into account, define $\vc{X}(t)$ as
\begin{align*}
  \vc{X}(t) = RC \Big(\Lambda^{\frac 12} \,\diag(\{\sigma_{sk}\}) \vc{B}_{s}^{*}(t) & + M Q \big(\diag(\{\lambda_{ek}^{\frac 32} \sigma_{ek}\}) \vc{B}_{e}^{*}(t)\\
  & \qquad + \sum_{k' \in K} \lambda_{ak'}^{\frac 12} \Gamma_{k'} \vc{B}_{k'}^{*}(t) \big)\Big) + t R \vc{c}.
\end{align*}
Since $R^{(n)} \equiv (I+G^{(n)})^{-1}$ is well defined for sufficiently large $n$, and for such $n$, \eqn{dif W 2} is can be written as
\begin{align}
\label{eqn:dif W 3}
  \widehat{\vc{W}}^{(n)}(t) = & \vc{X}^{(n)}(t) + R^{(n)} \widehat{\vc{Y}}^{(n)}(t),
\end{align}
where
\begin{align*}
  \vc{X}^{(n)}(t) =&  R^{(n)} C \big(\widehat{\vc{V}}^{(n)}_{s} (\ol{\vc{N}}^{(n)}_{a}(t)) + M^{(n)} Q (\widehat{\vc{N}}^{(n)}_{e}(t) + \sum_{k' \in K} \widehat{\vc{\Phi}}^{(n)}_{k'}(\ol{N}^{(n)}_{dk'}(t)) \big) \nonumber\\
  & \hspace{15ex} - R^{(n)} CM^{(n)} Q P^{\rs{t}} \widehat{\vc{\epsilon}}^{(n)}(t) + \sqrt{n} t R^{(n)} (\vc{\rho}^{(n)} - \vc{1}).
\end{align*}

It can be shown that $\vc{X}^{(n)}(\cdot) \overset{w}{\longrightarrow} \vc{X}(\cdot)$. Furthermore, $\widehat{Y}^{(n)}_{k}(t)$ increases only when $\widehat{W}^{(n)}_{k}(t) = 0$. Similar to the reflection mapping, if $R$ is a completely $\sr{S}$-matrix, this suggests that \eqn{dif W 3} implies that there are unique nonnegative $\vc{W}(\cdot)$ and $\vc{Y}(\cdot)$ such that
\begin{align}
\label{eqn:dif W 4}
  \vc{W}(t) = & \vc{X}(t) + R \vc{Y}(t),
\end{align}
This implication is actually proved in \cite{Will1998b}. One also needs to check the tightness for $\vc{W}^{(n)}(\cdot) \overset{w}{\longrightarrow} \vc{W}(\cdot)$ and $\vc{Y}^{(n)}(\cdot) \overset{w}{\longrightarrow} \vc{Y}(\cdot)$.

\begin{theorem}[Theorem 7.1 of \cite{Will1998a}] {\rm
\label{thr:multi-class GJN}
Under the assumptions \eqn{heavy 1} and \eqn{mul variance 1}--\eqn{dif Phi k 1}, if the service discipline is FCFS and if $I + G$ is invertible and its inverse is a completely $\sr{S}$-matrix, then we have
\begin{align*}
  (\vc{X}^{(n)}(\cdot), \vc{W}^{(n)}(\cdot), \vc{Y}^{(n)}(\cdot), \vc{L}^{(n)}(\cdot)) \overset{w}{\longrightarrow} (\vc{X}(\cdot), \vc{W}(\cdot), \vc{Y}(\cdot), \vc{L}(\cdot)).
\end{align*}
Furthermore, this convergence still holds for any head of line service discipline as long as the state collapse condition holds for some sequence $\Delta^{(n)}$ converging $\Delta$.
}\end{theorem}

For the multi-class GJN, there is no results about the continuity of the stationary distribution. The difficulty for this case is remarked in Section 5 of \cite{GamaZeev2006}.

%------------------------------------------------------ 
\section{Halfin-Whitt regime}
\label{sect:Halfin-Whitt}
%------------------------------------------------------ 

As we have observed in \sectn{extension to many servers}, the diffusion approximation by SRBM is still possible for the GJN with many server nodes. However, it ignores the influence of the number of idle servers which varies from 0 to the total number of servers at one node. So, the approximation may be too coarse. On the other hand, if there are infinitely many servers at a node, one may apply the central limit theorem to get asymptotics for approximation. In this case, there is no queue, and this is another extreme. Furthermore, the limiting process is Gaussian, and therefore not a diffusion process. So, we will not discuss about them. Comprehensive discussions about such a limiting process is found in Section 9 of \cite{Glyn1990}.

One may wonder what approximation is possible between those two extremes. This motivated Halfin and Whitt \cite{HalfWhit1981} to increase the number $n$ of servers as well as the arrival rate of customers in a single queue with many servers while keeping service time distributions. In this case, a network model may be too complicate because different nodes may have different numbers of customers. So, this type of limits has been studied mainly for a single node, that is, a single queue. We here refer to results in \cite{HalfWhit1981} for such a model.

Assume that the $n$-th system has $n$ servers, and let $L^{(n)}(t)$ be the number of customers in system at time $t$. Customers arrives according to a renewal process, and their service times are $i.i.d.$. An arriving customer has to wait if all servers are busy. An important assumption here is that the service time distribution is exponential. This is the crucial assumption which makes analysis possible and results tractable.

We will use the following notations. $\lambda^{(n)}$ is the arrival rate of customers, $\mu$ is the service rate, $\rho^{(n)} = \lambda^{(n)}/(n\mu)$. We assume the following heavy traffic condition. There is a constant $\beta$ such that
\begin{align}
\label{eqn:heavy 2}
  \lim_{n \to \infty} \sqrt{n} (1 - \rho^{(n)}) = \beta, \qquad \beta > 0.
\end{align}
We further assume the following condition on the moments of the inter-arrival time $T^{(n)}(\ell)$ for the $n$-th system.
\begin{align}
\label{eqn:variance 1}
 &\lim_{n \to \infty} (\lambda^{(n)})^{2} Var(T^{(n)}(1)) = c^{2},\\
\label{eqn:third 1}
 & \sup_{n \ge 1} (\lambda^{(n)})^{3} \dd{E}((T^{(n)}(1))^{3}) < \infty.
\end{align}
Define the following diffusion scaling for the number of customers in system.
\begin{align*}
  \widetilde{L}^{(n)}(t) = \frac 1{\sqrt{n}} (L^{(n)}(t) - n).
\end{align*}

The following theorems are obtained by Halfin and Whitt \cite{HalfWhit1981}.

\begin{theorem}[Theorem 3 of \cite{HalfWhit1981}]
\label{thr:HW 1}
  Assume \eqn{heavy 2}--\eqn{third 1}. If $\widetilde{L}^{(n)}(0) \overset{w}{\longrightarrow} Z(0)$, then $\widetilde{L}^{(n)}(\cdot) \overset{w}{\longrightarrow} Z(\cdot)$, where $Z(\cdot) \equiv \{Z(t); t \ge 0\}$ is a diffusion process characterized by the following stochastic integration.
\begin{align}
\label{eqn:diffusion process 1}
  Z(t) = Z(0) + \int_{0}^{t} \sigma(Z(u)) dB(u) + \int_{0}^{t} m(Z(u)) du, \qquad t \ge 0,
\end{align}
where $B(\cdot)$ is the standard Brownian motion, and $\sigma^{2}(x)$ and $m(x)$ are given by
\begin{align*}
  m(x) = - \mu (\beta + \mu x 1(x < 0)), \qquad \sigma^{2}(x) = \mu(1+c^{2}).
\end{align*}
\end{theorem}

\begin{theorem}[Theorem 4 of \cite{HalfWhit1981}]
\label{thr:HW 2}
  Assume \eqn{heavy 2}--\eqn{third 1}, and let $\widetilde{L}^{(n)}(\infty)$ and $Z(\infty)$ be random variables subject to the stationary distributions of $\widetilde{L}^{(n)}(\cdot)$ and $Z(\cdot)$, respectively, of \thr{HW 1}, then $\widetilde{L}^{(n)}(\infty) \overset{w}{\longrightarrow} Z(\infty)$. Furthermore, let $\beta^{*} = 2\beta/(1+c^{2})$ and let $\Phi_{N}(x)$ be the standard normal distribution function, then
\begin{align*}
  \lim_{n \to \infty} \dd{P}(L^{(n)} \ge n) = \dd{P}(Z(\infty) \ge 0) = (1 + \sqrt{2\pi} \beta^{*} \Phi_{N}(\beta^{*}) exp((\beta^{*})^{2}/2))^{-1},
\end{align*}
and 
\begin{align*}
 & \dd{P}(Z(\infty) > x|Z(\infty)\ge 0) = e^{-x \beta^{*}}, \qquad x > 0,\\
 & \dd{P}(Z(\infty) \le x|Z(\infty) \le 0) = \frac{\Phi_{N}(\beta^{*}+x)} {\Phi_{N}(\beta^{*})}, \qquad x \le 0.
\end{align*}
\end{theorem}

This approximation is popular for the performance evaluation of call centers, and there have been several attempts for its generalization concerning the service time distribution (e.g., see \cite{GamaGold2013,Reed2009}). In particular, Reed \cite{Reed2009} identifies the limiting process $Z(\cdot)$ for $\widetilde{L}^{(n)}(\cdot)$ as the solution of a certain integral equation. However, it is hard to see the property of this solution. So, Gamarnik and Goldberg \cite{GamaGold2013,Reed2009}) derive bounds for the tail probability of the stationary distribution.

%------------------------------------------------------ 
\section{Concluding Remarks}
\label{sect:concluding}
%------------------------------------------------------ 

From the diffusion approximations in Sections \sect{single-class}, \sect{diffusion} and \sect{multi-class}, we can see the following facts.
\begin{mylist}{2}
\item [(\sect{concluding}a)] For both single-class and multi-class queueing networks, SRBM is obtained as a process limit of the sequence of diffusion scaling processes.
\item [(\sect{concluding}b)] The SRBM for a single-class GJN is determined by the asymptotic means and variances of the counting processes $N_{ei}(\cdot)$ and $N_{si}(\cdot)$ and random routing function $\Phi_{ij}(\cdot)$.
\item [(\sect{concluding}c)] The SRBM for a multi-class GJN requires the asymptotic mean and variance of the accumulated workload $V_{si}(\cdot)$ in addition to those for a single-class queueing network.
\end{mylist}

For (\sect{concluding}b), we can get those data of $N_{ei}(\cdot)$ and $N_{si}(\cdot)$ from their logarithmic moment generating functions as defined in \eqn{gamma e}. Those logarithmic generating functions play important roles in the stationary equation of the GJN, which is the reason why they arise in the tail decay rate of the stationary distribution in \sectn{continuity tail}.

Although SRBM is very important for diffusion approximation, very little is known about its stationary distribution for the more than two dimensional case. Thus, it is often to use a product form stationary distribution, which is obtained only when the so called skew symmetric condition holds. However, a choice of the primitive data is very limited under this condition. Thus, it would be very nice if we can find some features of the stationary distribution.

The author and his colleagues \cite{DaiMiyaWu2015} recently studied to decompose the stationary distribution of SRBM into two marginals. They show that, if the stationary distribution is decomposable, then the marginal distributions are obtained as lower dimensional SRBM, and their primitive data are identified. This suggests that we may have more flexible choice of the primitive data than the skew symmetric case, so we may have a better approximation. 

There are two interesting observations in \sectn{extension to queue length} and \app{decay}. One is the fact that time scaling is not essential to get the limit of a sequence of the stationary distributions for diffusion compatible approximation. This may be useful for simulation to get the stationary distribution of the SRBM. Another is that large deviations rate functions are closely related to diffusion approximation through logarithmic moment generating function. Both observations may not be surprising, but it is notable that they are mathematically confirmed. 

\section*{Acknowledgements}
  I am grateful to Hideaki Yamashita of Tokyo Metropolitan University for inviting me to write this paper. I appreciate invaluable comments and suggestions from Ken'ichi Kawanishi of Gunma University. I also thank anonymous referees and Masahiro Kobayashi of Science University of Tokyo for their thoughtful comments. This work is supported by Japan Society for the Promotion of Science under grant No.\ 24310115.

\newpage
\appendix

\section*{Appendix}

\medskip

\section{Proof of \lem{Taylor 1}}
\label{app:Taylor 1}
We only prove the first formula of \eqn{Taylor 1} because \eqn{Taylor 2} can be obtained similarly. By the assumption \eqn{1-node condition 0}, $\widetilde{F}_{e}(w)$ with complex variable $w$ is analytic in a neighborhood of the origin in the complex plain $\dd{C}$. Since
\begin{align*}
  \frac {\partial }{\partial w} \big( e^{z} \widetilde{F}_{e}(w) - 1\big) \big|_{(z,w) = (0,0)} = \widetilde{F}'_{e}(0) = \dd{E}(T_{e}) \ne 0,
\end{align*}
$e^{z} \widetilde{F}_{e}(w) - 1=0$ has a unique root $w = w(z)$ which is analytic in $z$ in a neighborhood of the origin by the implicit function theorem (e.g., see Theorem 3.11 of the Volume II of \cite{Mark1977}). Hence, $\gamma_{e}(z)$ is analytic in a neighborhood of the origin, and therefore, it follows from Taylor expansion that
\begin{align*}
  \gamma_{e}(\theta) = - \eta(\theta) = - \Big( \eta(0) + \eta'(0) \theta + \frac 12 \eta''(0) \theta^{2} + o(\theta^{2}) \Big), \qquad \theta \to 0.
\end{align*}
Obviously, $\eta(0) = 0$. By differentiating both sides of $e^{\theta} \widetilde{F}_{e}(\eta) = 1$, we have
\begin{align}
\label{eqn:Taylor a1}
  \widetilde{F}_{e}'(\eta) \frac {d \eta}{d\theta} + \widetilde{F}_{e}(\eta) = 0,
\end{align}
where $\eta$ is a function $\theta$, namely, $\eta(\theta)$, but we drop $\theta$ for simplicity. This implies that $\eta'(0) = - \lambda_{e}$. Differentiating \eqn{Taylor a1}, we have
\begin{align}
\label{eqn:Taylor a2}
  \widetilde{F}_{e}'(\eta) \frac {d^{2} \eta}{d\theta^{2}} + \widetilde{F}_{e}''(\eta) \Big(\frac {d \eta}{d\theta}\Big)^{2} + \widetilde{F}_{e}'(\eta) \frac {d \eta}{d\theta} = 0,
\end{align}
which yields that $\eta''(0) = - \lambda_{e}^{3} \sigma_{e}^{2}$. Thus, we obtain \eqn{Taylor 1}.

\section{Proof of \thr{1-node L asym}}
\label{app:theorem L asym}

Let $\beta^{(n)} = 1 - \rho^{(n)}$. Since \eqn{1-node extra 1} and \eqn{1-node extra 2} hold, \eqn{1-node SE 3} yields, for $\theta \le 0$, 
\begin{align}
\label{eqn:1-node SD 1}
  \lim_{n \to \infty} \varphi^{(n)}(\beta^{(n)} \theta, -\lambda^{(n)}_{e} \beta^{(n)} \theta + o(\beta^{(n)} \theta), \lambda^{(n)}_{s} \beta^{(n)} \theta + o(\beta^{(n)} \theta)) = \frac {2} {2 - \lambda_{e}^{2} (\sigma_{e}^{2} + \sigma_{s}^{2}) \theta},
\end{align}
where the condition \eqn{1-node asym c3} is tacitly used to the small order terms $o(\beta^{(n)} \theta)$ as $n \to \infty$. To complete the proof, we need to show that
\begin{align}
\label{eqn:1-node SD 2}
  \lim_{n \to \infty} \varphi^{(n)}(\beta^{(n)} \theta, 0,0) = \frac {2} {2 - \lambda_{e}^{2} (\sigma_{e}^{2} + \sigma_{s}^{2}) \theta}, \qquad \theta \le 0.
\end{align}
For this, denote the right-hand side of \eqn{1-node SD 2} by $f(\theta)$. Then, by a triangular inequalities, we have
\begin{align*}
  |\varphi^{(n)} & (\beta^{(n)} \theta, 0,0) - f(\theta)| \\
  & \le |\varphi^{(n)}(\beta^{(n)} \theta, -\lambda^{(n)}_{e} \beta^{(n)} \theta + o(\beta^{(n)} \theta), \lambda^{(n)}_{s} \beta^{(n)} \theta + o(\beta^{(n)} \theta)) - f(\theta)|\\
 & \quad + |\varphi^{(n)}(\beta^{(n)} \theta, 0,0) - \varphi^{(n)}(\beta^{(n)} \theta, -\lambda^{(n)}_{e} \beta^{(n)} \theta + o(\beta^{(n)} \theta), \lambda^{(n)}_{s} \beta^{(n)} \theta + o(\beta^{(n)} \theta))|. 
\end{align*}
From \eqn{1-node SD 1}, the first absolute term in the right-hand side of this inequality vanishes as $n \to \infty$. Thus, we only need to show that its second absolute term vanishes as $n \to \infty$. For this, we write down it in terms of expectations for $\theta \le 0$.
\begin{align*}
  \mbox{The second absolute term} & = \Big| \dd{E}\Big( e^{\beta^{(n)} \theta L} \Big(1 - e^{(-\lambda^{(n)}_{e} \beta^{(n)} \theta + o(\beta^{(n)} \theta)) R_{e} + (\lambda^{(n)}_{s} \beta^{(n)} \theta + o(\beta^{(n)} \theta)) R^{(n)}_{s}} \Big)\Big) \Big|\\
  & \le \Big| 1 - e^{(-\lambda^{(n)}_{e} \beta^{(n)} \theta + o(\beta^{(n)} \theta)) R_{e}} \Big| \\
  & \quad + \Big| e^{(-\lambda^{(n)}_{e} \beta^{(n)} \theta + o(\beta^{(n)} \theta)) R_{e}} \Big(1 - e^{(\lambda^{(n)}_{s} \beta^{(n)} \theta + o(\beta^{(n)} \theta)) R^{(n)}_{s}} \Big)\Big) \Big|.
\end{align*}
The right-hand side converges to $0$ by the bounded convergence theorem because $\theta \le 0$ and $\beta^{(n)} > 0$ converges to $0$. This completes the proof of \thr{1-node L asym} because the convergence of Laplace transfers implies the convergence of the corresponding distributions.

\section{Normalizing factor in \thr{counting}}
\label{app:Donsker counting}

We check the normalizing factor to be correct in \thr{counting}. For this, we derive a limiting distribution for the one-dimensional marginal. It follows from the definition of $N^{(n)}$ that
\begin{align*}
  N^{(n)}(nt) > [kt] \quad \mbox{if and only if} \quad \sum_{\ell=1}^{[kt]} \tau^{(n)}_{\ell} < nt.
\end{align*}
Hence,
\begin{align*}
  \frac 1{\sqrt{nt}} & (N^{(n)}(nt) - \lambda^{(n)} nt)  > \frac 1{\sqrt{nt}} ([kt] - \lambda^{(n)} nt) \quad \\
  & \mbox{if and only if} \quad \frac 1{\sqrt{[kt]}\sigma^{(n)}} \sum_{\ell=1}^{[kt]} (\tau^{(n)}_{\ell} - m^{(n)}) < \frac 1{\sqrt{[kt]} \sigma^{(n)}} (nt - m^{(n)}[kt]).
\end{align*}
We let $n, k \to \infty$ so that $\frac 1{\sqrt{[kt]} \sigma^{(n)}} (nt - m^{(n)}[kt])$ converges to an arbitrarily given $x$. Then, 
\begin{align*}
  \lambda^{(n)} nt - [kt] \sim \lambda^{(n)} \sqrt{[kt]} \sigma^{(n)} x,
\end{align*}
where we recall that $\lambda^{(n)} = 1/m^{(n)}$, and therefore
\begin{align*}
 & \frac 1{\sqrt{nt}} ([kt] - \lambda^{(n)} nt) \sim - \frac{\lambda^{(n)} \sqrt{[kt]} \sigma^{(n)} x} {\sqrt{(\sqrt{[kt]} \sigma^{(n)} x + m^{(n)}[kt])t}} \sim - \frac{\lambda^{(n)} \sigma^{(n)} x} {\sqrt{m^{(n)}}} = - \sqrt{(\lambda^{(n)})^{3}} \sigma^{(n)} x,\\
 & [kt] \sim \lambda^{(n)} nt \sim \lambda nt
\end{align*}
Thus, we have
\begin{align*}
  \lim_{n \to \infty} \dd{P}\Big( \frac 1{\sqrt{nt}} (N^{(n)}(nt) - & \lambda^{(n)} nt) > - \sqrt{(\lambda^{(n)})^{3}} \sigma^{(n)} x \Big) \\
  & = \lim_{n \to \infty} \dd{P}\Big( \frac 1{\sqrt{[\lambda nt]}\sigma^{(n)}} \sum_{\ell=1}^{[\lambda nt]} (\tau^{(n)}_{\ell} - m^{(n)}) < x \Big).
\end{align*}
Since the right-hand side is the standard normal distribution function, $\sigma \lambda^{\frac 32}$  correctly appears in the theorem.

\section{Proof of \thr{decay rate continuity}}
\label{app:decay}

We only need to prove that \eqn{gamma continuity} implies the convergence of $\sr{D}^{(n)}$ to $\widetilde{\sr{D}}$. A major problem is to show that $1/\sqrt{n}$ in the second formula of \eqn{gamma continuity} is a right scaling factor. For this, we first derive the stationary equation similar to \eqn{1-node SE 1} for a general GJN. The reason why we consider general $d$ rather than $d=2$ is that there are not so much difference between them and the general case may have own interest. For $i,j \in J$ and $\vc{\theta}, \vc{\eta}, \vc{\xi} \in \dd{R}^{d}$, let
\begin{align*}
 & \varphi(\vc{\theta}, \vc{\eta}, \vc{\xi}) = \dd{E}(e^{\br{\vc{\theta}, \vc{L}(0)} + \br{\vc{\eta}, \vc{R}_{e}(0)} + \br{\vc{\xi}, \vc{R}_{s}(0)}}), \\
 & \varphi_{i0}(\vc{\theta}, \vc{\eta}, \vc{\xi}) = \dd{E}(e^{\br{\vc{\theta}, \vc{L}(0)} + \br{\vc{\eta}, \vc{R}_{e}(0)} + \br{\vc{\xi}, \vc{R}_{s}(0)}}1(L_{i}(0)=0)), \\
 & \varphi_{ei}(\vc{\theta}, \vc{\eta}, \vc{\xi}) = \dd{E}_{ei}(e^{\br{\vc{\theta}, \vc{L}(0-)} + \br{\vc{\eta}, \vc{R}_{e}(0-)} + \br{\vc{\xi}, \vc{R}_{s}(0-)}}), \\
 & \varphi_{ei0-}(\vc{\theta}, \vc{\eta}, \vc{\xi}) = \dd{E}_{ei}(e^{\br{\vc{\theta}, \vc{L}(0-)} + \br{\vc{\eta}, \vc{R}_{e}(0-)} + \br{\vc{\xi}, \vc{R}_{s}(0-)}}1(L_{i}(0-)=0)),\\
 & \varphi_{si}(\vc{\theta}, \vc{\eta}, \vc{\xi}) = \dd{E}_{si}(e^{\br{\vc{\theta}, \vc{L}(0-)} + \br{\vc{\eta}, \vc{R}_{e}(0-)} + \br{\vc{\xi}, \vc{R}_{s}(0-)}}), \\
 & \varphi_{si0+}(\vc{\theta}, \vc{\eta}, \vc{\xi}) = \dd{E}_{si}(e^{\br{\vc{\theta}, \vc{L}(0+)} + \br{\vc{\eta}, \vc{R}_{e}(0+)} + \br{\vc{\xi}, \vc{R}_{s}(0+)}}1(L_{i}(0+)=0)),\\
 & \varphi_{sij0}(\vc{\theta}, \vc{\eta}, \vc{\xi}) = \dd{E}_{si}(e^{\br{\vc{\theta}, \vc{L}(0+)} + \br{\vc{\eta}, \vc{R}_{e}(0+)} + \br{\vc{\xi}, \vc{R}_{s}(0+)}}1(L_{i}(0+)=0, L_{j}(0-)=0)),
\end{align*}
where $\dd{E}_{ui}$ represents the expectation with respect to the Palm measure concerning $N_{ui}$ for $u=e,s$, which is defined similarly to \eqn{Palm 1}. We then apply the rate conservation law similarly to \eqn{RCL path} for the following $X(t)$.
\begin{align*}
  X(t) = \left\{\begin{array}{ll}
  g_{i}(\vc{\theta}) e^{\br{\vc{\theta}, \vc{L}(t)} + \br{\vc{\eta}, \vc{R}_{e}(t)} + \br{\vc{\xi}, \vc{R}_{s}(t)}}, \quad & L_{i}(t) = 0,\\
  e^{\br{\vc{\theta}, \vc{L}(t)} + \br{\vc{\eta}, \vc{R}_{e}(t)} + \br{\vc{\xi}, \vc{R}_{s}(t)}}, \quad & L_{i}(t) \ge 1.
  \end{array} \right.
\end{align*}
This yields the following formula for $\vc{\theta} \le \vc{0}$, where $\theta_{0} = 0$.
\begin{align}
\label{eqn:GJN SE 1}
   - \sum_{i \in J} & (\eta_{i} + \xi_{i}) \varphi(\vc{\theta}, \vc{\eta}, \vc{\xi}) + \sum_{i \in J} (\eta_{i} (1-g_{i}(\vc{\theta})) + \xi_{i}) \varphi_{i0}(\vc{\theta}, \vc{\eta}, \vc{\xi}) \nonumber\\
   & + \sum_{i \in J} \lambda_{ei} (e^{\theta_{i}} \widetilde{F}_{ei}(\eta_{i}) - 1) \varphi_{ei}(\vc{\theta}, \vc{\eta}, \vc{\xi}) \nonumber\\
   & + \sum_{i \in J} \lambda_{ei} \big(e^{\theta_{i}}  \widetilde{F}_{ei}(\eta_{i}) ( \widetilde{F}_{si}(\xi_{i})-1) + 1-g_{i}(\vc{\theta})) \varphi_{ei0-}(\vc{\theta}, \vc{\eta}, \vc{\xi} \big) \nonumber\\
  & + \sum_{i \in J} \lambda_{si} \Big(\sum_{j \in J \cup \{0\}} e^{-\theta_{i}+\theta_{j}} p_{ij} \widetilde{F}_{si}(\xi_{i}) - 1\Big) \varphi_{si}(\vc{\theta}, \vc{\eta}, \vc{\xi}) \nonumber \\
  & + \sum_{i \in J} \sum_{j \in J \cup \{0\}} \lambda_{si} e^{-\theta_{i}+\theta_{j}} p_{ij} (g_{i}(\vc{\theta}) - \widetilde{F}_{si}(\xi_{i}) ) \varphi_{si0+}(\vc{\theta}, \vc{\eta}, \vc{\xi}) \nonumber \\
  & + \sum_{i \in J} \sum_{j \in J} \lambda_{si} e^{-\theta_{i}+\theta_{j}} p_{ij} g_{i}(\vc{\theta}) (\widetilde{F}_{sj}(\xi_{j}) - 1) \varphi_{sij0}(\vc{\theta}, \vc{\eta}, \vc{\xi}) = 0.
\end{align}

We next choose $\eta_{i}, \xi_{i}$ and $g_{i}(\vc{\theta})$ in such a way that
\begin{align*}
  e^{\theta_{i}} \widetilde{F}_{ei}(\eta_{i}) = 1, \qquad \widetilde{F}_{si}(\xi_{i}) \sum_{j \in J \cup \{0\}} e^{-\theta_{i}+\theta_{j}} p_{ij} = 1, \qquad g_{i}(\vc{\theta}) = \widetilde{F}_{si}(\xi_{i}), \qquad i \in J.
\end{align*}
Denote these $\eta_{i}$ and $\xi_{i}$ by $\eta_{i}(\theta_{i})$ and $\xi_{i}(\vc{\theta})$, and let $\vc{\eta}(\vc{\theta})$ and $\vc{\xi}(\vc{\theta})$ be their vectors, then we arrive at
\begin{align}
\label{eqn:GJN SE 2}
   - \sum_{i \in J} & (\eta_{i}(\theta_{i}) + \xi_{i}(\vc{\theta})) \varphi(\vc{\theta}, \vc{\eta}(\vc{\theta}), \vc{\xi}(\vc{\theta})) \nonumber\\
   & + \sum_{i \in J} (\eta_{i}(\theta_{i}) (1-\widetilde{F}_{si}(\xi_{i}(\vc{\theta})))+ \xi_{i}(\vc{\theta})) \varphi_{i0}(\vc{\theta}, \vc{\eta}(\vc{\theta}), \vc{\xi}(\vc{\theta})) \nonumber\\
  & + \sum_{i , j \in J, i \ne j} \lambda_{si} e^{-\theta_{i}+\theta_{j}} p_{ij} \widetilde{F}_{si}(\xi_{i}(\vc{\theta})) (\widetilde{F}_{sj}(\xi_{j}(\vc{\theta})) - 1) \varphi_{sij0}(\vc{\theta}, \vc{\eta}(\vc{\theta}), \vc{\xi}(\vc{\theta})) = 0.
\end{align}
For $d=2$, it follows from \eqn{logarithmic 1} and \eqn{gamma d} that
\begin{align*}
  \eta_{i}(\theta_{i}) = - \gamma_{ei}(\theta_{i}), \qquad \xi_{i}(\vc{\theta}) = - \gamma_{di}(\vc{\theta}).
\end{align*}

We now consider the $n$-th GJN with $d=2$, and its characteristics are indexed by superscript ``$^{(n)}$''.
Hence, recalling the definition of $\gamma_{+}$, \eqn{GJN SE 2} becomes, for $\vc{\theta} \le \vc{0}$,
\begin{align}
\label{eqn:GJN SE 3}
   \gamma^{(n)}_{+}(\vc{\theta}) \varphi^{(n)}(&\vc{\theta}, - \vc{\gamma}^{(n)}_{e}(\vc{\theta}), - \vc{\gamma}_{d}^{(n)}(\vc{\theta})) - \sum_{i=1,2} \gamma^{(n)}_{di}(\vc{\theta}) \varphi^{(n)}_{i0}(\vc{\theta}, - \vc{\gamma}^{(n)}_{e}(\vc{\theta}), - \vc{\gamma}_{d}^{(n)}(\vc{\theta})) \nonumber\\
   & - \sum_{i = 1,2} \gamma^{(n)}_{ei}(\theta_{i}) (1-t_{i}^{-1}(\vc{\theta})) \varphi_{i0}^{(n)}(\vc{\theta}, - \vc{\gamma}^{(n)}_{e}(\vc{\theta}), - \vc{\gamma}_{d}^{(n)}(\vc{\theta})) \nonumber\\
  & + \sum_{(i,j) = (1,2), (2,1)} \lambda^{(n)}_{si} e^{-\theta_{i}+\theta_{j}} p_{ij} t_{i}^{-1}(\vc{\theta}) \nonumber\\
  & \hspace{10ex} \times (t_{j}^{-1}(\vc{\theta}) - 1) \varphi_{sij0}^{(n)}(\vc{0}, - \vc{\gamma}^{(n)}_{e}(\vc{\theta}), - \vc{\gamma}_{d}^{(n)}(\vc{\theta})) = 0,
\end{align}
where we have used the fact that
\begin{align*}
  \varphi_{sij0}^{(n)}(\vc{\theta}, - \vc{\gamma}^{(n)}_{e}(\vc{\theta}), - \vc{\gamma}_{d}^{(n)}(\vc{\theta})) = \varphi_{sij0}^{(n)}(\vc{0}, - \vc{\gamma}^{(n)}_{e}(\vc{\theta}), - \vc{\gamma}_{d}^{(n)}(\vc{\theta})),
\end{align*}
because $\vc{L}(0) = \vc{0}$ implies that $\br{\vc{\theta}, \vc{L}(0)} = 0 = \br{\vc{0}, \vc{L}(0)}$. Thus, using the diffusion scaling $\widehat{\gamma}^{(n)}_{ui}(\vc{\theta})$ of \eqn{logarithmic diffusion} and substituting $\frac 1{\sqrt{n}} \vc{\theta}$ into the first variable $\vc{\theta}$ of $\varphi^{(n)}, \varphi^{(n)}_{i0}, \varphi^{(n)}_{sij0}$ in \eqn{GJN SE 3}, we have, under diffusion scaling, 
\begin{align}
\label{eqn:GJN SE 4}
   \widehat{\gamma}^{(n)}_{+}&(\vc{\theta}) \varphi^{(n)} \Big( \frac 1{\sqrt{n}} \vc{\theta}, - \frac 1n \widehat{\vc{\gamma}}^{(n)}_{e}(\vc{\theta}), - \frac 1n \widehat{\vc{\gamma}}_{d}^{(n)}(\vc{\theta}) \Big) \nonumber\\
   & - \sum_{i=1,2} \widehat{\gamma}^{(n)}_{di}(\vc{\theta}) \varphi^{(n)}_{i0}\Big(\frac 1{\sqrt{n}}\vc{\theta}, - \frac 1n \widehat{\vc{\gamma}}^{(n)}_{e}(\vc{\theta}), - \frac 1n \widehat{\vc{\gamma}}_{d}^{(n)}(\vc{\theta})\Big) \nonumber\\
   & - \sum_{i = 1,2} \widehat{\gamma}^{(n)}_{ei}(\theta_{i}) \Big(1-t_{i}^{-1}\Big(\frac 1{\sqrt{n}} \vc{\theta}\Big)\Big) \varphi_{i0}^{(n)}\Big(\frac 1{\sqrt{n}} \vc{\theta}, - \frac 1n \widehat{\vc{\gamma}}^{(n)}_{e}(\vc{\theta}), - \frac 1n \widehat{\vc{\gamma}}_{d}^{(n)}(\vc{\theta})\Big) \nonumber\\
  & + \sum_{(i,j) = (1,2), (2,1)} n \lambda^{(n)}_{si} e^{-\theta_{i}+\theta_{j}} p_{ij} t_{i}^{-1}\Big(\frac 1{\sqrt{n}} \vc{\theta}\Big) \nonumber\\
  & \hspace{12ex} \times \Big(t_{j}^{-1}\Big(\frac 1{\sqrt{n}} \vc{\theta}\Big) - 1\Big) \varphi_{sij0}^{(n)}\Big(\vc{0}, - \frac 1n \widehat{\vc{\gamma}}^{(n)}_{e}(\vc{\theta}), - \frac 1n \widehat{\vc{\gamma}}_{d}^{(n)}(\vc{\theta})\Big) = 0.
\end{align}

Since $\frac 1n \widehat{\vc{\gamma}}_{ui}^{(n)}(\vc{\theta}) \to 0$ as $n \to \infty$ for $u = e,d$ by \eqn{gamma continuity}, we have, using similar arguments in \app{theorem L asym},
\begin{align*}
 & \lim_{n \to \infty} \varphi^{(n)} \Big( \frac 1{\sqrt{n}} \vc{\theta}, - \frac 1n \widehat{\vc{\gamma}}^{(n)}_{e}(\vc{\theta}), - \frac 1n \widehat{\vc{\gamma}}_{d}^{(n)}(\vc{\theta}) \Big) = \tilde{\varphi}(\vc{\theta}).
\end{align*}
Since $t_{j}^{-1}\Big(\frac 1{\sqrt{n}} \vc{\theta}\Big) - 1$ has the order $1/\sqrt{n}$ as $n \to \infty$, this and the fact that \eqn{GJN SE 4} converges to the BAR \eqn{BAR 2} of the limiting SRBM concludes
\begin{align*}
 & \lim_{n \to \infty} \lambda^{(n)}_{si} \sqrt{n} \varphi^{(n)}_{i0} \Big( \frac 1{\sqrt{n}} \vc{\theta}, - \frac 1n \widehat{\vc{\gamma}}^{(n)}_{e}(\vc{\theta}), - \frac 1n \widehat{\vc{\gamma}}_{d}^{(n)}(\vc{\theta}) \Big) = \tilde{\varphi}_{i}(\vc{\theta}),
\end{align*}
and therefore $-(\lambda^{(n)}_{si} \sqrt{n})^{-1} \widehat{\vc{\gamma}}_{di}^{(n)}(\vc{\theta})$ is a right term to converge to $\tilde{\gamma}_{i}(\vc{\theta})$ so that the continuity of the stationary distribution holds. Since $-\widehat{\gamma}_{di}^{(n)}(\vc{\theta})$ and $-\gamma_{di}^{(n)}(\vc{\theta})$ have the same sign and 
\begin{align*}
  -\gamma_{di}^{(n)}(\vc{\theta}) = (\tilde{F}^{(n)}_{si})^{-1}(t^{-1}_{i}(\vc{\theta})),
\end{align*}
$-\gamma_{di}^{(n)}(\vc{\theta}) < 0$ if and only if $t_{i}(\vc{\theta}) > 1$, which is the condition to determine $\Gamma_{3-i}$ for \thr{decay rate 1}.  
This completes the proof of \thr{decay rate continuity}.

%------ corresponding author ------------------ 
\vspace{15pt}
{\leftskip 8.5cm \parindent 0mm
Masakiyo Miyazawa\\
Department of Information Sciences\\
Tokyo University of Science\\
Yamazaki 2641, Noda\\
Chiba 278-8510, Japan\\ 
E-mail: \texttt{miyazawa@rs.tus.ac.jp}
\par}
%------------------------------------------------------ 
\end{document}